\documentclass[pdflatex,sn-mathphys-num]{sn-jnl}
\usepackage{graphicx}%
\usepackage{multirow}%
\usepackage{amsmath,amssymb,amsfonts}%
\usepackage{amsthm}%
\usepackage{mathrsfs}%
\usepackage[title]{appendix}%
\usepackage{xcolor}%
\usepackage{textcomp}%
\usepackage{manyfoot}%
\usepackage{booktabs}%
\usepackage{algorithm}%
\usepackage{algorithmicx}%
\usepackage{algpseudocode}%
\usepackage{listings}%
\usepackage{subfigure}
\usepackage{enumitem}
\newtheorem{lemma}{Lemma}
 
\newtheorem{theorem}{Theorem}%  meant for continuous numbers
\newtheorem{definition}{Definition}
\newcommand{\rev}[1]{{#1}}
\newcommand{\ReLU}{\operatorname{ReLU}}
\newcommand{\clconv}{\operatorname{cl\,conv}}
\newtheorem{remark}{Remark}%
\usepackage{makecell} % 导言区
\usepackage[utf8]{inputenc}
\usepackage[T1]{fontenc}

\usepackage{geometry} % 引入宏包
\begin{document}

\title[An Efficient Stochastic Subgradient Method for Global Placement in VLSI Circuits]{An Efficient Stochastic \rev{Subgradient} Method for the Global Placement Problem in \rev{Very Large-Scale Integration Circuits}}
%%=============================================================%%
%% GivenName	-> \fnm{Joergen W.}
%% Particle	-> \spfx{van der} -> surname prefix
%% FamilyName	-> \sur{Ploeg}
%% Suffix	-> \sfx{IV}
%% \author*[1,2]{\fnm{Joergen W.} \spfx{van der} \sur{Ploeg} 
%%  \sfx{IV}}\email{iauthor@gmail.com}
%%=============================================================%%

\author[2,1]{\fnm{Yi-Shuang} \sur{Yue}}\email{yueyishuang@amss.ac.cn}
\author[1,2]{\fnm{Yu-Hong} \sur{Dai}}\email{dyh@lsec.cc.ac.cn}
%\equalcont{These authors contributed equally to this work.}
\author*[1,2]{\fnm{Haijun} \sur{Yu}}\email{hyu@lsec.cc.ac.cn}
%\equalcont{These authors contributed equally to this work.}

\affil[1]{\orgdiv{State Key Laboratory of Mathematical Sciences (SKLMS) and State Key Laboratory of Scientific and Engineering Computing (LSEC), Institute of Computational Mathematics and Scientific/Engineering Computing}, \orgname{Academy of Mathematics and Systems Science, Chinese Academy of Sciences}, \orgaddress{\street{No. 55 Zhongguancun East Road}, \city{Haidian District}, \postcode{100190}, \state{Beijing}, \country{China}}}

\affil[2]{\orgdiv{School of Mathematical Sciences}, \orgname{University of Chinese Academy of Sciences}, \orgaddress{\city{Beijing}, \postcode{100049}, \country{China}}}

%\affil[2]{\orgdiv{Department}, \orgname{Organization}, \orgaddress{\street{Street}, \city{City}, \postcode{10587}, \state{State}, \country{Country}}}

%\affil[3]{\orgdiv{Department}, \orgname{Organization}, \orgaddress{\street{Street}, \city{City}, \postcode{610101}, \state{State}, \country{Country}}}

%%==================================%%
%% Sample for unstructured abstract %%
%%==================================%%
\abstract{The placement problem in Very Large-Scale Integration (VLSI) circuits is a critical step in chip design. Its primary goal is to optimize the wirelength of circuit components within a confined area while adhering to nonoverlapping constraints. Most analytical placement models rely on smooth approximations, thereby sacrificing the accuracy of wirelength estimation. To mitigate these inaccuracies, this paper introduces a novel approach that directly optimizes the original nonsmooth wirelength and proposes an innovative penalty model tailored for the global placement problem. Specifically, we transform the nonoverlapping constraints into rectified linear penalty functions, allowing for a more precise formulation of the problem. Notably, we recast the resultant optimization problem into a form analogous to training deep neural network with Rectified Linear Units (ReLU). Leveraging automatic differentiation techniques from deep learning, we efficiently compute the subgradient of the objective function. This facilitates the application of stochastic subgradient methods to solve the model. To enhance the algorithm's performance, several advanced techniques are further introduced, leading to significant improvements in both efficiency and solution quality. Numerical experiments were conducted on \rev{Gigascale Systems Research Center (GSRC)} benchmark and \rev{International Symposium on Physical Design 2005 (ISPD2005)} benchmark circuits. The results demonstrate that our proposed model and algorithm achieve significant reductions in wirelength while effectively eliminating overlaps. This highlights the potential of our approach as a transformative advancement for VLSI placement.
\rev{Furthermore, we establish a rigorous convergence proof for the proposed stochastic subgradient method. To the best of our knowledge,  it constitutes the first such result for the ReLU-type nonsmooth and nonconvex optimization problems.}}


%%\pacs[JEL Classification]{D8, H51}

%%\pacs[MSC Classification]{35A01, 65L10, 65L12, 65L20, 65L70}

\maketitle

\section{Introduction}\label{sec1}
Circuit placement is a fundamental and critical stage in the physical design of very large-scale integration (VLSI) circuits. Its primary goal is to determine the specific locations of physical units within a predefined placement area, ensuring that device nonoverlapping constraints are satisfied while optimizing key performance metrics such as wirelength, delay, and area utilization~\cite{alpert2008handbook}. The quality of circuit placement has a profound impact on subsequent design stages, including post-placement optimization and routing, directly influencing circuit performance, power consumption, chip area, manufacturing cost, signal integrity, and so on. The placement problem is a typical NP-hard problem~\cite{kahng2011vlsi}, necessitating the development of advanced algorithms to address its complexity. 

Existing solutions can broadly be classified into three categories: simulated annealing~\cite{sechen2012vlsi}, partitioning~\cite{lasalle2013multi}, and analytical methods~\cite{viswanathan2007fastplace}. The simulated annealing method emulates the physical annealing process and gradually optimizes the placement by randomly swapping, rotating, or moving modules, assessing the effects of these changes to find the global optimum. Partitioning methods address the placement problem by recursively dividing large instances—comprising net lists and placement regions into smaller sub-problems. Once the instance is sufficiently reduced, efficient local optimization algorithms are applied to achieve high-quality placements within each partition~\cite{caldwell1999optimal}. Analytical methods leverage linear or nonlinear programming techniques to obtain optimal or near-optimal solutions for the placement problem. By modeling placement objectives and constraints mathematically, analytical approaches have emerged as a dominant paradigm for solving VLSI placement problems, offering exceptional accuracy and scalability~\cite{huang2021optimization}.

The placement problem can be rigorously formulated as a constrained optimization problem, where the primary objective is to minimize the half-perimeter wirelength (HPWL) while satisfying nonoverlapping constraints for physical units. Different methodologies for approximating and optimizing wirelength have led to the emergence of two major categories of placement models: smooth models and nonsmooth models. The smooth model approximates the original nonsmooth HPWL using continuously differentiable objective functions, such as $l_p$ norms~\cite{chan2005multilevel}, log-sum-exp functions~\cite{naylor2001non}, and a weighted average model~\cite{hsu2011tsv}. To address the overlap constraints, smoothing techniques like the Gaussian equation ~\cite{chen2008ntuplace3} and the Helmholtz equation~\cite{chan2005multilevel} are often employed. These techniques relax the problem by transforming the original nonsmooth optimization into a smooth nonlinear programming problem. Smooth models are widely adopted in modern placement solvers, which incorporate penalty methods, conjugate gradient methods, and Nesterov's accelerated techniques for efficient optimization.  The nonsmooth optimization models directly address the original nonsmooth HPWL. The HPWL objective function inherently involves absolute values, which introduce nonsmooth and nondifferentiable points, making it a challenging optimization landscape. To overcome these difficulties, researchers have developed specialized algorithms tailored for nonsmooth structures, such as the conjugate subgradient method~\cite{zhu2015nonsmooth} and mixed-variable optimization~\cite{sun2023floorplanning}, which are more adept at handling discontinuities in the gradient caused by absolute values.

Existing smooth models, while exhibiting high computational efficiency, suffer from limitations in accuracy due to their reliance on smooth approximations. This limitation becomes particularly pronounced in large-scale designs, where approximation errors can accumulate and lead to suboptimal solutions and additional adjustment costs. On the other hand, nonsmooth models can accurately represent the original structure of the HPWL, but their nondifferentiable objective functions and constraints make conventional gradient-based optimization methods inapplicable. Therefore, the core motivation of our research is to design a solution algorithm that preserves the accuracy of HPWL representation while handling nondifferentiable points and maintaining high computational efficiency. 
In recent years, a growing line of work has explored enhancing global placement with 
high-performance computing and data-driven methods. GPU-based frameworks such as DREAMPlace~3.0~\cite{gu2020dreamplace} and Xplace~\cite{liu2023xplace} exploit parallel modeling and neural components to achieve scalable and high-quality placement. More recently, transferable learning has been introduced through graph neural networks, as in TransPlace~\cite{hou2025transplace}, which leverages placement knowledge across designs. Beyond these, large language models 
(LLMs) have been applied to automatically discover optimization algorithms and heuristics~\cite{yao2025evolution,sun2025automatically}, pointing to an emerging paradigm of AI-driven automation in placement research.

In this paper, we propose a promising nonsmooth optimization model for VLSI global placement. In this model, the objective function is defined as the HPWL wirelength function, and the nonoverlapping constraint is characterized by the piecewise-linear hat function related to the overlap of the elements, which is free of gradient vanishing and gradient exploding if gradient methods are used. Due to the nonconvex and nonsmooth properties of the problem, we design a stochastic subgradient (specifically, the generalized gradient of Clarke) descent operator splitting algorithm to solve it. The solution procedure is formulated as a neural net training problem and is implemented using the deep-learning-based toolkit PyTorch~\cite{paszke2019pytorch}. Several dedicated techniques are further introduced, including degree-weighted sampling based on hypergraphs, adaptive parameter updating, and random disturbance to improve the performance and placement effect of the algorithm. In addition, we propose an alternative-optimization legalization algorithm by combining wirelength optimization and the de-overlap technique, which eliminates overlap while maintaining the advantage of a short wirelength. Moreover, we provide a theoretical convergence analysis of the proposed stochastic subgradient algorithm, establishing its reliability for nonsmooth optimization in large-scale placement problems. We conducted numerical experiments to compare our proposed algorithm with the benchmark provided by the Gigascale Systems Research Center (GSRC)~\cite{Gsrc}. The results demonstrate significant improvements in wirelength reduction and overlap elimination. In addition, large-scale experiments on the International Symposium on Physical Design 2005 Placement Benchmarks (ISPD2005) show that our method achieves wirelength and overlap performance comparable to state-of-the-art placement tools. %although it needs more optimization due to the lack of GPU acceleration and explicit density optimization.  
These results highlight the advantages and potential of our proposed model and algorithms in solving VLSI placement problems.

The rest of this article is organized as follows. In the second section, we first present the mathematical model for the placement problem, followed by the derivation of the corresponding penalty function model. The third section is devoted to the stochastic subgradient method to optimize the nonsmooth penalty model. The fourth section provides a convergence analysis, and the fifth section presents experimental results on the GSRC benchmark and the ISPD2005 benchmark.  Finally, the last section offers a summary of the article.

\section{The placement problem and mathematical models}\label{sec2}

\subsection{The problem definition}
A VLSI circuit board can be modeled as a hypergraph $G(V, E)$, where $V=\{v_1,v_2,\ldots,v_n\}$ denotes the set of cells, each cell $v_i$ represents a physical component of the circuit, and $E=\left\{e_{1}, \ldots,e_{m}\right\}$ denotes the set of nets. Each net $e_k \in E$ is a subset of $V$, corresponding to a group of cells that are electrically interconnected.

Each cell $v_i \in V$, for $i = 1, \ldots, n$, is modeled as a rectangular block defined by its width $w_i$ and height $h_i$.  Let $(x_i, y_i)$ denote the center coordinates of cell $v_i$. The placement region for these cells is a designated rectangular area bounded by the coordinates $(0,0)$ at the lower left corner and $(W, H)$ at the upper right corner.

The objective of the VLSI placement problem is to determine the optimal positions for the cells within the placement region such that the total wirelength, measured as the sum of the lengths of the nets connecting the cells, is minimized. Additionally, the placement must satisfy the constraint that no two cells overlap within the designated area.
A mathematical formulation of the VLSI global placement problem for mixed cells is described as follows~\cite{chang2009essential}:
\begin{align}
      \min _{z} \mathrm{HPWL}(z)& :=\sum_{e \in E}\Bigl( \max _{\left(v_{i}, v_{j}\right) \in e}\left|x_{i}-x_{j}\right|+\max _{\left(v_{i}, v_{j}\right) \in e}\left|y_{i}-y_{j}\right| \Bigr) , \label{eq:HPWL} \\
      \text { s.t. } \quad &\left|x_{i}-x_{j}\right| \geq \frac{w_{i}+w_{j}}{2} \text { or }\left|y_{i}-y_{j}\right| \geq \frac{h_{i}+h_{j}}{2}, \quad \forall v_{i}, v_{j} \in V, \label{eq:overlap} \\
      & \frac{w_{i}}{2} \leq x_{i} \leq W-\frac{w_{i}}{2}\ \text{ and } \frac{h_{i}}{2} \leq y_{i} \leq H-\frac{h_{i}}{2}, \quad \forall(x_{i}, y_{i}) \in V. 
      \label{eq:bound}
\end{align}
Here $z=(x, y)$, $x=(x_1, \ldots, x_n)$, $y=(y_1, \ldots, y_n)$.
HPWL stands for half-perimeter wirelength. The objective function and constraints specified in the problem (\ref{eq:HPWL})--\eqref{eq:bound} are nonconvex and nonsmooth, making it a challenging optimization problem. To overcome these difficulties, a penalty method is employed to relax the constraints, enabling the reformulation and transformation of the problem into a neural network framework. The reformulated problem is then optimized using advanced deep-learning-based techniques.

\subsection{A rectified linear penalty method}

To enforce the boundary condition, we adopt the following penalty terms:
\begin{equation}\label{eq:penalty_b}
\begin{aligned}   
    \ell_{b}\left(x_{i}, y_{i}\right)&:=\operatorname{ReLU}\Bigl(\frac{w_{i}}{2}-x_{i}\Bigr)+\operatorname{ReLU}\left(x_{i}-w_{i}^{+}\right) \\
    &\qquad +\operatorname{ReLU}\Bigl(\frac{h_{i}}{2}-y_{i}\Bigr)+\operatorname{ReLU}\left(y_{i}-h_{i}^{+}\right),
\end{aligned}
\end{equation}
where $w_{i}^{+}=W-w_{i} / 2$, $h_{i}^{+}=H-h_{i} / 2$, $\operatorname{ReLU}(x):=\max (0, x)$.
To enforce the nonoverlapping constraints, a straightforward approach introduces the following penalty function:
\begin{gather}
    \ell_{x}\left(x_{i}, x_{j}\right):=\operatorname{ReLU}\left(\frac{w_{i}+w_{j}}{2}+x_{i}-x_{j}\right) \times\operatorname{ReLU}\left(\frac{w_{i}+w_{j}}{2}+x_{j}-x_{i}\right), \label{eq:lx} \\
    \ell_{y}\left(y_{i}, y_{j}\right):=\operatorname{ReLU}\left(\frac{h_{i}+h_{j}}{2}+y_{i}-y_{j}\right) \times\operatorname{ReLU}\left(\frac{h_{i}+h_{j}}{2}+y_{j}-y_{i}\right). \label{eq:ly}
\end{gather}
This leads to the following penalty formulation:
\begin{equation}
    \min _{z}\biggl\{\operatorname{HPWL}(z)+\sum_{\left(x_{i}, y_{i}\right) \in V} \gamma_i\ell_{b}\left(x_{i}, y_{i}\right)+\sum_{i<j}\gamma_{ij} \ell_{x}\left(x_{i}, x_{j}\right) \ell_{y}\left(y_{i}, y_{j}\right)\biggr\}, \label{eq:model0}
\end{equation}
where $\gamma_i$ is the penalty parameter of term $\ell_{b}\left(x_{i}, y_{i}\right)$, and $\gamma_{ij}$ is the penalty parameter of term $\ell_{x}\left(x_{i}, x_{j}\right) \ell_{y}\left(y_{i}, y_{j}\right)$. However, due to the peculiarities of overlapping penalty terms (the generalized directional derivative of $\ell_x \ell_y$ can be very close to $0$), this penalty function model may encounter the problem of gradient vanishing during optimization. 

Alternative approaches exist for constructing the overlap penalty. In~\cite{zhu2015nonsmooth}, the placement area is uniformly partitioned into a grid of nonoverlapping bins, and the squared difference between the overlapping cell area and the bin area is used as the penalty term. A recent work~\cite{sun2023floorplanning}, on the other hand, chooses the square root of the overlapping area as the penalty term. However, both of these methods suffer from the gradient vanishing issue, which limits their effectiveness in optimization.

To remove the deficiency of the nonoverlapping penalty, we need to use different penalties for $\ell_{x}\ell_{y}$ in \eqref{eq:model0}. We note that (\ref{eq:lx}) and (\ref{eq:ly}) are quadratic polynomials when (\ref{eq:overlap}) are not satisfied, which implies that the derivatives of the objective can be very small. To remedy this, we propose the following two-dimensional piecewise-linear hat function to replace the piecewise-quadratic function $\ell_x\ell_y$:
\begin{equation} \label{eq:2dhatfun}
    \varphi_{r,t}(x, y) =
    \begin{cases}
        0,                 & \text{if } |x| \ge r \text{ or } |y| \ge t; \\[4pt]
        1 - \dfrac{|x|}{r}, & \text{if } \dfrac{|y|}{t} \leq \dfrac{|x|}{r} < 1; \\[4pt]
        1 - \dfrac{|y|}{t}, & \text{if } \dfrac{|x|}{r} < \dfrac{|y|}{t} < 1.
    \end{cases}
\end{equation}
Then, a two-dimensional penalty for pairwise overlap can be constructed as
\begin{equation}\label{eq:new-penalty}
    \ell_{x y}^{}\left(x_{i}, y_{i}, x_{j}, y_{j}\right)=\varphi_{w_{i j}, h_{i j}}\left(x_{i}-x_{j}, y_{i}-y_{j}\right), \quad w_{i j}=\frac{w_{i}+w_{j}}{2}, \ h_{i j}=\frac{h_{i}+h_{j}}{2}.
\end{equation}
The new penalty problem can  be modeled as
\begin{equation}
    \min _{z} f(z):=\biggl\{\operatorname{HPWL}(z)+ \sum_{\left(x_{i}, y_{i}\right) \in V} \gamma_i \ell_{b}\left(x_{i}, y_{i}\right)+ \sum_{i<j} \gamma_{ij} \ell_{x y}^{}\left(x_{i}, y_{i}, x_{j}, y_{j}\right)\biggr\}.
    \label{eq:model1}
\end{equation}

\rev{\begin{remark}
Note that the third term in \eqref{eq:model1} involves pairwise nonoverlapping penalties, which may result in a massive number of terms for large-scale designs. In practice, these terms are not evaluated exhaustively at each iteration. Instead, we adopt a stochastic, batch-based optimization strategy (detailed in Section \ref{degree-Weighted Sampling}) to select a subset of critical nets and their associated interactions. This degree-weighted sampling strategy significantly reduces computational costs while preserving solution quality, making the proposed approach scalable to large-scale placement problems.
\end{remark}
}
\begin{remark}
We can also replace $1-|x|/r$ with $(rt-t|x|)$ and $1-|y|/t$ with $(rt-r|y|)$ in \eqref{eq:2dhatfun}, but our numerical experiments show that they lead to similar results.
\end{remark}

\begin{figure}[ht]
\centering
\includegraphics[width=0.8\textwidth]{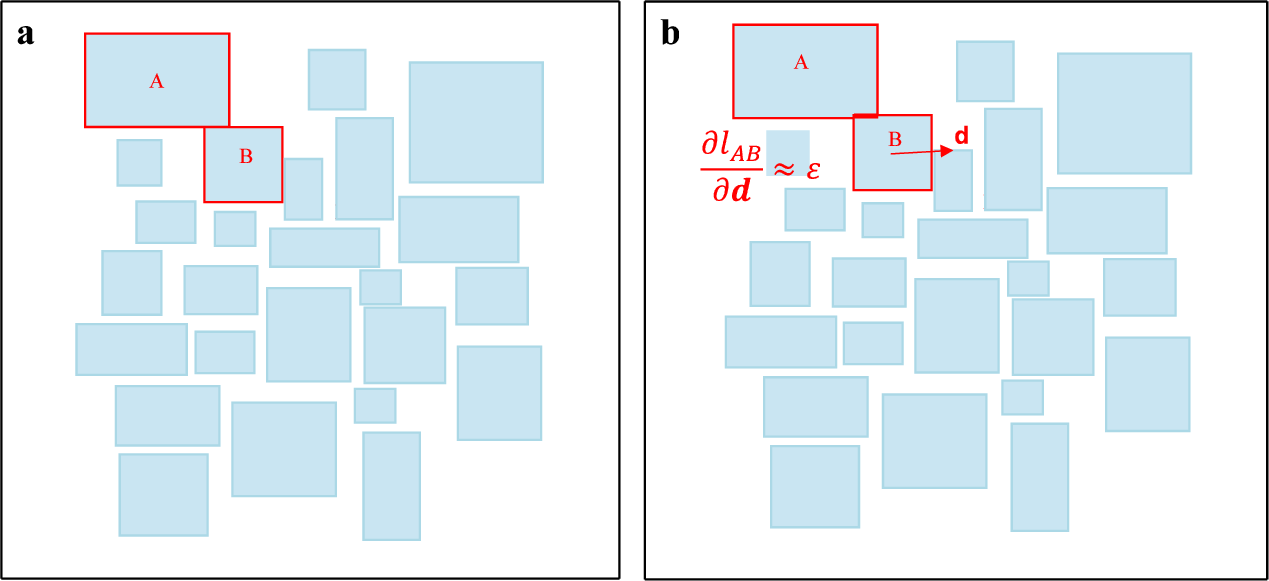}
\caption{{\bf a)} A placement corresponding to a local minimum point of the original problem \eqref{eq:HPWL}-\eqref{eq:bound}; {\bf b)} a small perturbation of the placement {\bf a}, where along certain directions, the HPWL decreases finitely, while the directional derivative of the overlap penalty may be arbitrarily small. Therefore, the objective function  of the penalty model \eqref{eq:model1} will decrease along this direction, and the nonoverlapping constraint will be violated, no matter how large the penalty constant is. In this case, the local minimum of model \eqref{eq:model1} will not be a local minimum of the original problem.}\label{bound}
\end{figure}

\bigskip
We note that it can be proven that, with a sufficiently large penalty constant, a one-dimensional rectified linear penalty is exact. However, the two-dimensional rectified linear penalty for pairwise overlap is not. Let us consider a case presented in Figure~\ref{bound}. In the left plot, the derivatives with respect to the positions of $A$ and $B$ cells are zero.
However, in the right plot, a directional derivative of the pairwise nonoverlapping penalty
with respect to the position of $B$ can be infinitely close to 0; e.g., the directional derivative along a direction $d$ which is a small perturbation to the $x$-direction. Therefore, if the HPWL has a finite nonzero gradient along this direction, the nonoverlapping constraint will be violated. Thus, the penalty is not exact.

Although the new penalty function model is not exact, we can prove that the amplitude of nonzero Clarke subgradients of the penalty term has both a lower and an upper bound (see Theorem \ref{gradnorm-bound}), thus avoiding the problems of gradient vanishing and gradient exploding, and better optimization results can be obtained.
To show this, we first prove that the new penalty model \eqref{eq:model1} can be transformed into a deep ReLU network.

\begin{theorem} \label{thm:equiv-to-relu-net}
The new model $f(z)$ defined in \eqref{eq:model1} can be formulated as a deep neural network with $\max\{\lceil \log M\rceil + 1, 3\}$ ReLU layers, where $M = \max_{e\in E} \{\#e\}$
and $\# e$ denote the number of elements in a net $e$.
\end{theorem}

\begin{proof}
First, for the absolute value function, we can reformulate it as
$$|x|=\operatorname{ReLU}(x)+\operatorname{ReLU}(-x).$$
The maximum of two values can be expressed as a ReLU net with one-hidden layer as
\begin{equation}
    \label{eq:abs-relu}
    \max(a,b)=\operatorname{ReLU}(a-b)+b = \operatorname{ReLU}(b-a)+a.
\end{equation}
The maximum of 4 values can be expressed as a composition of two layers of ReLU functions by \eqref{eq:abs-relu} and the following identity
$$\max(a, b, c, d) = \max(\max(a,b),\max(c,d)).$$
 Similarly, the maximum of $2^n$ values can be expressed as a composition of $n$ layers of ReLU functions. Thus, the HPWL defined in \eqref{eq:HPWL} can be formulated into a deep ReLU neural network with $\lceil \log M \rceil+1$ ReLU layers. Here $\lceil x \rceil$ represents the smallest integer that is no smaller than $x$.
 The boundary penalty \eqref{eq:penalty_b} is already formulated as a network of one ReLU layer. The pairwise penalty \eqref{eq:2dhatfun} can be written as 
 $$ \varphi_{r,t}(x, y) =\max\bigl\{ \text{ReLU}(1 -{|x|}/{r}), \text{ReLU}(1 - {|y|}/{t}) \bigr\}.$$
 Thus, it can be formulated as a neural network with 3 ReLU layers. 
 The summations in \eqref{eq:model1} only affect the width of the network without influencing its depth.
 The theorem is proved.
 \end{proof}

\medskip

Next, we define the subgradient used in this paper.
\begin{definition}[Clarke subgradient~\cite{clarke1990optimization}]
For a locally Lipschitz function $F$, the Clarke subdifferential is
\[
\partial^{C}F(z)=\clconv\Big\{\lim_{j\to\infty}\nabla F(z_j):\ 
z_j\to z,\ F \text{ differentiable at } z_j\Big\},
\]
which, for a finite piecewise-affine $F$, equals the convex hull of the gradients
of the active pieces at $z$.
An arbitrarily chosen element in $\partial^C F(z)$ is referred to as a (Clarke) subgradient of $F$ at $z$.
\end{definition}
Note that for convex and Lipschitz functions, the Clarke subgradient is equivalent to the usual subgradient. For differentiable functions, the Clarke subgradient is equivalent to the usual gradient.  
In this sense, we may omit the superscript $C$ in $\partial^C F(z)$ for notation simplicity.

\begin{theorem} \label{gradnorm-bound}

    Define $S=\{z:\ 0\in\partial f(z)\}$, where $\partial f(z)$ denotes the Clarke subdifferential of the objective defined in model \eqref{eq:model1}.
    Then 
    \begin{equation}\label{eq:non-zero-gradient-set}
    G_{S^c} := \{\, \|s\| \mid s \in \partial f(z),\, z\notin S \,\}
    \end{equation}
    has a positive lower and upper bound.
   % From the definition of Clarke subgradient we can know the $l^p$ norm of nonzero gradients of the penalty term $\sum_{\left(x_{i}, y_{i}\right) \in V} \gamma_i \ell_{b}\left(x_{i}, y_{i}\right)+ \sum_{i<j} \gamma_{ij} \ell_{x y}^{}\left(x_{i}, y_{i}, x_{j}, y_{j}\right)$ has a lower bound and an upper bound.
   
\end{theorem} 
\begin{proof}

% We consider the $l^1$ norm. The cases $p\neq 1$ are similar. 
% The $l^1$ norm of the nonzero subgradient of $ \ell_{b}\left(x_{i}, y_{i}\right):=\operatorname{ReLU}\left(\frac{w_{i}}{2}-x_{i}\right)+\operatorname{ReLU}\left(x_{i}-w_{i}^{+}\right)+\operatorname{ReLU}\left(\frac{h_{i}}{2}-y_{i}\right)+\operatorname{ReLU}\left(y_{i}-h_{i}^{+}\right)$  satisfies
% $$ 1 \le \|\partial\ell_{b}(x,y)\|_1 \le 2.
% $$
% The $l^1$ norm of the nonzero subgradient of $l_{xy}(x_i,y_i,x_j.y_j)$
% defined in \eqref{eq:new-penalty} is:
% $$ \bigl\|\partial\ell_{xy}(x_i,y_i,x_j.y_j) \bigr\|_1= \begin{cases}
%   \displaystyle \frac{1}{w_{ij}}, &\displaystyle \frac{|y_i-y_j|}{h_{ij}} \leq \frac{|x_i-x_j|}{w_{ij}}< 1,\\
%     \displaystyle \frac{1}{h_{ij}}, &\displaystyle \frac{|x_i-x_j|}{w_{ij}}< \frac{|y_i-y_j|}{h_{ij}} < 1.
% \end{cases}
% $$
% Thus it is clear that the nonzero gradient norm of the penalty term $\sum_{\left(x_{i}, y_{i}\right) \in V} \gamma_i \ell_{b}\left(x_{i}, y_{i}\right)+ \sum_{i<j} \gamma_{ij} \ell_{x y}^{}\left(x_{i}, y_{i}, x_{j}, y_{j}\right)$ has a lower and an upper bound.

i) Since $f(z)$ is continuous and piecewise-affine with a finite number of pieces, it is Lipschitz continuous. Consequently, $G_{S^c}$ is bounded above by $L_g$, defined as the maximum norm of the gradient vectors across all affine pieces.

ii) Since there is only a finite number of distinct gradients, the norm of the gradient at any differentiable point $z\notin S$ has a uniform positive lower bound. Since the Clarke subdifferential $\partial f(z)$ at a given point $z$ is a convex hull of a finite number of given gradients, if $0\notin \partial f(z)$, then by the separating hyperplane argument~\cite{boyd2004convex}, there exists a unit vector $u_z$ and $\beta_z>0$, such that
\begin{equation}
\langle v,u_z\rangle\ \ge\ \beta_z \qquad  \forall\, v\in\partial f(z).\label{sephyparg}
\end{equation}
Since there are only a finite number of affine regions with nonzero gradients to define subdifferentials, we can choose a finite number of different values $\beta_z$ to make \eqref{sephyparg} hold for any $z\notin S$. 
Then, there exists $\beta=\min \beta_z > 0$. Thus, we have $\| v\|> \beta$, for any $v\in \partial f(z)$, $z\notin S$,
i.e., $G_{S^c}$ has a lower bound $\beta$. The theorem is proved.

\end{proof}

\smallskip

\subsection{Analogy to deep learning}
\begin{figure}[ht]
\centering
\includegraphics[width=1\textwidth]{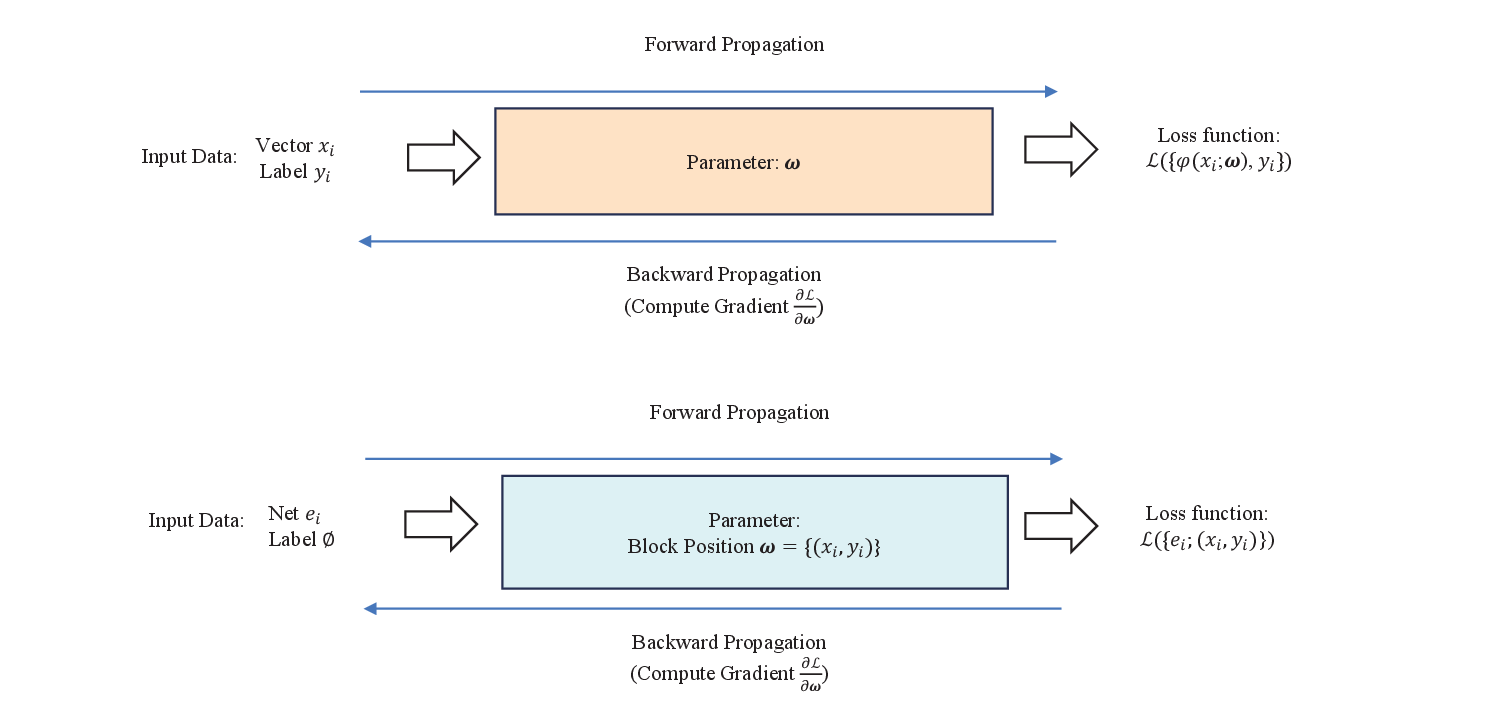}
\caption{Analogy between neural network training and placement problem training. The top plot (a) illustrates a supervised learning process, and the bottom plot (b) illustrates our method, which is a reinforcement learning approach. Both of them utilize back propagation to solve the network training problem. }\label{Analogy}
\end{figure}

Both analytic placement and training neural nets are essentially solving nonlinear optimization problems. In the context of supervised learning, neural network training involves feeding input data, consisting of feature vectors and corresponding labels, into the network~\cite{lin2019dreamplace}. The network computes a loss function based on the predicted and actual values, and optimization algorithms are employed to iteratively update the network parameters using gradients. This process continues until the loss is minimized and the network accurately predicts the labels for the given data.
  In our placement problem, we directly optimize the objective function composed of wirelength and overlap boundary constraint penalties,
   continuously updating the block positions until the objective function is optimized, which can be viewed as a reinforcement learning approach.    
  Furthermore, a direct analogy can be drawn between the placement problem and neural network training; see Fig.~\ref{Analogy}. The input data instance is replaced with net information; the feature vector represents key information such as net connection and block information. The label is set to an empty set, and the block positions $\{(x_i, y_i)\}$ correspond to net parameters. The loss is defined as our objective functions (\ref{eq:model1}), and the parameters, i.e., the block positions, are adjusted throughout the training process to minimize the loss function.  By establishing this one-to-one mapping, the placement problem is effectively reformulated as a neural net training problem. This transformation allows the use of deep-learning-based techniques to solve placement problems, using forward propagation to compute objectives and backward propagation to calculate gradients, thereby enabling the implementation of efficient optimization strategies.

\section{The stochastic algorithm}\label{sec3}

\subsection{Random batch splitting method}
In the last section, the original constrained optimization problem was reformulated into an unconstrained optimization problem using penalty methods. Due to the nonsmooth and nonconvex nature of the problem, traditional optimization techniques may struggle to efficiently identify solutions close to the global minimum.
Therefore, we will employ a stochastic subgradient method to solve this nonsmooth optimization problem. 
Moreover, we will employ the operator splitting method~\cite{jin2022random} to optimize the interaction function in placement problems. This method improves efficiency and precision by decomposing the objective function into two components: 
\begin{equation}
    f(z) = f_1(z)+f_2(z), \label{function}
\end{equation}
where
$$f_1(z)=\sum_{\left(x_{i}, y_{i}\right) \in V} \gamma_i \ell_{b}\left(x_{i}, y_{i}\right)
    + \sum_{i<j} \gamma_{ij} \ell_{x y}^{}\left(x_{i}, y_{i}, x_{j}, y_{j}\right)$$ 
    is a short-range component, and
    $$f_2(z)=\sum_{e \in E} \ell_{e}, \quad \ell_e:=\max _{\left(v_{i}, v_{j}\right) \in e}\left|x_{i}-x_{j}\right|+\max _{\left(v_{i}, v_{j}\right) \in e}\left|y_{i}-y_{j}\right| $$
    is a long-range component.
The short-range component $f_1$ mainly addresses the loss caused by boundaries and overlap. Since this term is nonzero only when a unit overlaps with others or exceeds the boundary, the calculation involving this component can be completed quickly. The computational complexity can be reduced from $O (N^2)$ to $ O (N)$ by filtering items that overlap with the block through meshing. The remote component $f_2$ is the long-range loss, the loss function is bounded, and the calculation amount can be effectively reduced using mini-batch techniques. 

In the subsequent sections, we provide a detailed explanation of the algorithm's implementation, focusing on how additional techniques are combined to effectively optimize the placement problem.

\subsubsection{Adaptive parameter adjustment}

In penalty-based optimization, the penalty parameters play a critical role in balancing the trade-off between the primary objective (e.g., wirelength minimization) and the constraint satisfaction (e.g., overlap avoidance). 
If the penalty parameters are too small, the algorithm may ignore constraint violations, producing inaccurate placements. Conversely, excessively large parameters can dominate the loss function, causing numerical instability and making the optimization trajectory overly sensitive to local perturbations, which may hinder convergence. Therefore, both the \emph{role} and the \emph{tuning} of penalty parameters are crucial for achieving stable and accurate solutions.

To address this issue, we adopt an \emph{adaptive} penalty parameter strategy inspired by~\cite{zhu2015nonsmooth,lu2014eplace}. The key idea is to dynamically increase the penalty for components that contribute most to the constraint violation while keeping other penalties moderate, so that the optimization focuses on the currently most critical violations without overwhelming the objective function.

Specifically, we initialize all penalty parameters with a uniform value $\gamma_0$. 
Then, at each iteration, the penalty parameters are updated as follows:
\begin{equation}
    \bar{\gamma_i} = 
    \begin{cases}
         \left\lceil \max \left(
            \frac{\partial (f_2)_i}{\partial (\ell_b(x_i, y_i))_i }
         \right) \right\rceil, & \text{if } 0 \notin \partial (\ell_b(x_i, y_i))_i,\\[3pt]
        \gamma_0, & \text{otherwise},
    \end{cases}
    \label{eq:gamma1}
\end{equation}
\begin{equation}
    \bar{\gamma_{ij}} = 
    \begin{cases}
         \left\lceil \max \left(
            \frac{\partial (f_2)_i}{\partial (\ell_{xy}(x_i, y_i, x_j, y_j))_i }
         \right) \right\rceil, & \text{if } 0  \notin \partial (\ell_{xy}(x_i, y_i, x_j, y_j))_i, \\[3pt]
        \gamma_0, & \text{otherwise}.
    \end{cases}
    \label{eq:gamma2}
\end{equation}
Here, $\partial(f_2)_i$ denotes the generalized derivatives of the wirelength term $f_2$ with respect to $(x_i,y_i)$, and similar definitions apply to the block term $\ell_{b}$ and the pairwise term $\ell_{xy}$. The element-wise division reflects the relative sensitivity of the objective to each constraint; $\max$ denotes taking the maximum over all directions. The ceiling function $\lceil \cdot \rceil$ ensures that the updated penalty parameter is nondecreasing in practice, which stabilizes the iteration and accelerates convergence.

This adaptive strategy effectively emphasizes the most critical violations at each step, allowing the algorithm to progressively reduce overlap and improve placement quality without requiring manual tuning of penalty schedules. 
Experimental results demonstrate that the use of rounding via the ceiling function significantly accelerates the algorithm, enhancing both computational efficiency and solution quality.

\subsubsection{Degree-weighted sampling based on a hypergraph}\label{degree-Weighted Sampling}
 
 To improve the efficiency further, we employ a degree-weighted sampling method based on a hypergraph to train the net in batches.  In this method, each net is represented as a vertex in the hypergraph, and an edge is established between two vertices if the corresponding nets share a common block.

Using this topological structure, the ``degree'' of each net vertex is defined as the number of connections it has with other vertices, reflecting its connectivity density~\cite{zhou2007bipartite}.  Higher-degree nets are more closely related to other nets, and their optimization is typically more challenging. To prioritize the optimization of critical nets, we adopt a sampling strategy that gives higher selection probability to nets with larger degrees. Specifically, for the net set $E = \{e_1, \ldots, e_m\}$, if the degree of the net $e_i$ is $d_i$, the sampled probability $p_i$  of net $e_i$ is calculated as follows: 
\begin{equation}
p_i=\frac{e^{d_i/T}}{\sum_{j=1}^m e^{d_j/T}}, \label{eq:probability}
\end{equation}
where $T$ is a constant parameter resembling temperature.
 This weighted sampling strategy not only improves computational efficiency but also ensures that sufficient attention is directed toward key nets. Consequently, it accelerates overall convergence while maintaining high solution quality.
The algorithm of degree-weighted sampling is as follows (Algorithm \ref{algorithm batch}):
\begin{algorithm}[ht]
\small{
    \caption{Degree-weighted sampling}\label{algorithm batch}
    \begin{algorithmic}[1]
        \State Create a mapping for each terminal to record which nets contain the terminal.
        \State Build the hypergraph of mapping and calculate the degree of the vertices.
        \State Calculate the sampling probability for each net using formula \eqref{eq:probability}.
        \State Sample from the full net $E$ with probability obtained from step 3.
    \end{algorithmic}
    }
\end{algorithm}

\subsubsection{Mean-field force}
To achieve a more compact arrangement of blocks, a regularization term known as ``mean-field force''~\cite{weiss1907hypothese} is incorporated into the objective function. This term is designed to reduce the distance between each block and the center of the placement region, thereby promoting tighter block placement. The regularization term is defined as:
$$
\mathcal{E}_{mean}=\sum_{v_i\in V}{\alpha \bigl\|\left( x_i,y_i \right) -\left( x_{mean},y_{mean} \right) \bigr\|^2_2}.
$$
Here, $\alpha$ is the force parameter. $\left( x_{mean},y_{mean} \right)$ represents the average position of all blocks. By incorporating this regularization term into the optimization process, the arrangement of blocks becomes more balanced, thereby reducing the overall placement area and significantly decreasing the total wirelength. Moreover, the inclusion of this term accelerates convergence to a stable configuration and mitigates the risk of being trapped in local minima.

\subsubsection{Incorporation of random disturbance}
One of the main challenges of nonconvex optimization problems is the potential to become trapped in local minima. To address this issue and promote a more comprehensive exploration of the solution space,
we have integrated a stochastic process into our algorithm. These stochastic components are designed to introduce randomness, thereby enhancing the algorithm's ability to escape local minima. The key stochastic components are implemented as follows:
\begin{enumerate}
    \item \textbf{Random initialization}: Initial coordinates $(x_i^0, y_i^0)$ are randomly sampled from a specified range to enhance the diversity of initial configurations.
    \item \textbf{Incorporation of random forcing}: By adding Gaussian random disturbances to the gradient, noise is introduced into the optimization process of the algorithm, which can effectively help the optimization escape local minima~\cite{xiao2021artificial}. The algorithm for adding random perturbations can be summarized as Algorithm \ref{algorithm disturbance}:
          \begin{algorithm}[ht]
          \small{
              \caption{Gaussian random disturbance}\label{algorithm disturbance}
              \begin{algorithmic}[1]
                  \For{Step $k$}
                  \State Calculate the subgradient $\boldsymbol{g}$ of loss function.
                  \State Generate a Gaussian random vector $\boldsymbol{\eta}$ with the same shape as the subgradient.
                  \State Calculate the perturbation coefficient by $\epsilon=0.2/k^3$.
                  \State Add  perturbation  and update subgradient by:
                  $$ \boldsymbol{g} = \boldsymbol{g} + \epsilon \|\boldsymbol{g}\| \boldsymbol{\eta}.$$
                  \EndFor
              \end{algorithmic}
              }
          \end{algorithm}
\end{enumerate}
By incorporating these stochastic ingredients, the algorithm has enhanced its efficiency in exploring the solution space, thereby increasing the likelihood of identifying globally optimal or near-optimal solutions.

\subsection{The overall random batch splitting algorithm}
Combining all the above details, our algorithm can be summarized as Algorithm~\ref{RBSM}, which we refer to as the random batch splitting method (RBSM).
\noindent
\begin{algorithm}[ht]
\small{
    \caption{RBSM for function \eqref{function}}\label{RBSM}
    \begin{algorithmic}[1]
       \State \textbf{Input:} $iter_{max}$, tolerance $\epsilon_{hpwl} = 0.0001$, $\epsilon_{overlap} = 2\%$;
        \State \textbf{Output:} Optimal placement position $(x^*, y^*)$;
        \State $k \gets 0$.
        \State Initialize $\bigl({x}^0, {y}^0\bigr)$ randomly and set initial learning rate $lr=lr_0 = 0.1$.
        \While{$k < iter_{max}$}
            \State Split $f(z) $ into two parts: distance term $f_1$ and penalty term $f_2$. 
            \State Calculate $\gamma_i^k$ and $\gamma_{ij}^k$ using equations \eqref{eq:gamma1}, \eqref{eq:gamma2}.
            \For{$step = 1$ to $25$}
            \State Sample from $ E$ by Algorithm \ref{algorithm batch} to obtain batch $E_q$. 
            % in the early stage, while  the full set $E$ in the later stage.
                \For{batch $E_q$}
                        \State Compute distance loss $f_1$ and get subgradient via automatic differentiation. 
						\State Add random perturbations by Algorithm \ref{algorithm disturbance} and update positions using a gradient method: 
						     $$({x},{y})\leftarrow(x,y)+lr*( g+\epsilon\| g\| {\mathbf \eta}).$$
						\State Compute penalty term $f_2$ and get subgradient via automatic differentiation.
						\State Add random perturbations to update positions using a gradient method similarly.     
                \EndFor
            \EndFor
            \State Get $\bigl(x^{k+1}, y^{k+1}\bigr) $ after all steps.
            \State $lr \gets lr_0 \times \left(1 + \cos\left(\frac{\pi \times k}{iter_{max}}\right)\right) / 2$.
            \If {$\left|\frac{\operatorname{HPWL}(x^{k}, y^{k})}{\operatorname{HPWL}(x^{k+1}, y^{k+1})} - 1\right| < \varepsilon_{\text{hpwl}}$  \textbf{and}  $\text{ratio}_{\text{overlap}} < \varepsilon_{\text{overlap}}$ }
                \State \textbf{Stop} and \textbf{return} $(x^*, y^*)=\bigl(x^k, y^k\bigr)$ .
            \EndIf
            \State $k \gets k + 1$
        \EndWhile
    \end{algorithmic}
    }
\end{algorithm}

In the position update phase, distance loss and penalty are calculated separately, enabling two targeted position updates. This approach allows for more effective and focused adjustments to overlap and wirelength.  After each update of $\gamma$, the model is trained until the local minimum is reached. Based on experimental results, it was observed that stochastic gradient descent (SGD) typically converges to a local minimum within 25 iterations. Therefore, a fixed iteration count of 25 is used instead of a variable stopping criterion. 
%Here, we use a fixed training step instead of a nonfixed step to the local minimum. This is because the experiment found that no more than 25 steps gradient descent method can reach the local minimum. 
To update the learning rate, we adopt the cosine annealing schedule from PyTorch, which gradually reduces the learning rate over iterations. This scheduling strategy helps improve convergence stability and avoid oscillations near local minima. At last, we  adopted a practical and reasonable termination criterion commonly used in analytical placement methods,  which is to stop the algorithm when the target HPWL satisfies condition $\left|\frac{\text{HPWL}(x^k, y^k)}{\text{HPWL}(x^{k+1}, y^{k+1})} - 1\right| < \varepsilon_{\text{hpwl}}$ and the overlap ratio meets  $\text{ratio}_{\text{overlap}} < \varepsilon_{\text{overlap}} = 2\%$.       
 
 The dominant cost is the penalty term $f_2$. With mini-batching, let $M$ be the
total number of nets, $b$ the batch size (nets per batch), $S=25$ the number
of inner steps per outer iteration, and $K$ the number of outer iterations.
Each inner step evaluates $f_2$ on one batch, costing $O(b^2)$ (pairwise
interactions inside the batch), while $f_1$ and perturbations cost $O(b)$.
Hence, the overall complexity is
\[
O\big(K S b^2\big).
\]
If the $S$ inner steps in one outer iteration visit disjoint batches that
together cover all nets (so $S=M/b$), this becomes
\[
O\big(K S b^2\big)=O\big(K M b\big),
\]
which is linear in $M$ for fixed $b$. Therefore, the mini-batched
implementation removes the quadratic dependence on the total number of nets,
and $f_2$ remains the primary computational bottleneck within each batch.

\subsection{Alternate-optimization legalization}
Due to the fact that the penalty is not exact, overlaps cannot be entirely eliminated during the global placement stage. Thus, further de-overlapping is required in the legalization phase. Typically, traditional legalization focuses solely on eliminating overlaps, which might lead to a significant increase in wirelength.  
 To address this, we use a strategy that alternates between optimizing overlap and minimizing distance loss. This method simultaneously reduces overlap and prevents significant increases in wirelength.  Initially, the algorithm calculates distance loss and updates positions to reduce the distance between components. It then assesses the overlap penalty, calculates gradients via automatic differentiation, and updates positions based on calculated steps to minimize overlaps. 
 This process of alternating updates continues until overlaps are eliminated and distance stabilization is achieved. 
 We note that for the rectified linear penalty, the legalization of nonoverlapping can be made exact in one step for individual pairwise constraints if these constraints are not coupled. However, due to the existence of coupled non-overlapping constraints in dense layouts, multiple rounds of legalization iterations are necessary.
 
 The specific steps are given in Algorithm \ref{algorithm overlap}.
\begin{algorithm}[!hpb]
\small{
    \caption{Legalization Algorithm}\label{algorithm overlap}
    \begin{algorithmic}[1]
        \For{$k$ from $1$ to $10$}
        \State Calculate distance loss by $\ell_{e}:=\max _{(v_{i}^k, v_{j}^k) \in e}\bigl|x_{i}^k-x_{j}^k\bigr|+\max _{(v_{i}^k, v_{j}^k) \in e}\bigl|y_{i}^k-y_{j}^k\bigr| .$
        \State Update positions by gradient descent method.
            \For{block $I(\text{position}(x_{i}, y_{i}))$ in $V$}
                \If{block $J(\text{position}(x_{j}, y_{j}))$ overlaps $I$}
                \State Get subgradient $g_i$ of $\ell_{x y}$ for $(x_{i}, y_{i})$  by automatic differentiation.
                \State Calculate step by\[
                    \alpha_i = \begin{cases}
                        r_{ij} - |x_i - x_j| & \text{if } \frac{\partial \ell_{x y}}{\partial x} \neq 0, \\
                        t_{ij} - |y_i - y_j| & \text{else if } \frac{\partial \ell_{x y}}{\partial y} \neq 0, \\
                        0 &\text{otherwise}.
                    \end{cases}
                \]
                %If multiple conditions hold simultaneously, one of them may be selected arbitrarily.
                
                \State update position of $(x_{i}, y_{i})$ by $(x_{i}^{k+1}, y_{i}^{k+1}) = (x_{i}, y_{i}) - \alpha_i *g_i/\|g_i\|$.
                \State \textbf{break}
                \EndIf
            \EndFor
            \State Calculate boundary loss $\ell_{b}\left(x_{i}, y_{i}\right)$ and get subgradient $g_i$ of $\ell_{b}(x_{i}, y_{i})$  by automatic differentiation.
            \State Calculate step by \[
                \alpha_i = \begin{cases}
                    x_i-(W - \frac{w_i}{2}) & \text{if } \frac{\partial \ell_b}{\partial x} > 0, \\
                    \frac{w_i}{2}-x_i     & \text{else if } \frac{\partial \ell_b}{\partial x} < 0, \\
                    y_i-(H - \frac{h_i}{2}) & \text{else if } \frac{\partial \ell_b}{\partial y} > 0, \\
                    \frac{h_i}{2}-y_i     & \text{else if } \frac{\partial \ell_b}{\partial y} < 0,\\
                    0 & \text{otherwise}.
                \end{cases}
            \]
            %If multiple conditions hold simultaneously, one of them may be selected arbitrarily.
            \State update position of $(x_{i}, y_{i})$ by $(x_{i}^{k+1}, y_{i}^{k+1}) = (x_{i}, y_{i}) - {\alpha}_i * g_i$.
        \EndFor
    \end{algorithmic}
    }
\end{algorithm}

\section{Convergence analysis of the RBSM algorithm}

We analyze the stochastic subgradient iteration underlying Algorithm \ref{RBSM}.
Let $f(z)$ be our nonsmooth, nonconvex objective with
$z=(x_1,y_1,\ldots,x_n,y_n)\in\mathbb{R}^{2n}$.
One outer iteration performs $M$ inner updates (we use $M=25$ in practice) on mini-batches $E_q$.
Index the overall run by $k=tM+m$ ($t\in\mathbb{N}$, $m\in\{0,\dots,M-1\}$) and set
$z_{k}=z_{t,m}$.
Each inner update  can be written as
\begin{equation}\label{eq:rbsm-one-step}
z_{k+1}=z_k-\alpha_k(g_k+\epsilon_k\|g_k\|\eta_k)=z_k-\alpha_k(g_k+\xi_k)=z_k-\alpha_k\hat{g}_k,
\end{equation}
where $g_k$ is a (mini-batch, importance-weighted) subgradient estimator of $f$ at $z_k$,
$\eta_k\sim\mathcal{N}(0,I)$ is independent of the past, and $\epsilon_k=0.2/k^3$ is the
perturbation coefficient from Algorithm \ref{algorithm disturbance}. 

Here we decompose our objective function into the form of the sum of boundary constraints and other functions, that is,
 $f=f_0+f_b$, where $f_b$ is finite boundary constraints 
\[
f_b(z)=\sum_{i=1}^n \gamma_i \big[ \operatorname{ReLU}\Bigl(\frac{w_{i}}{2}-x_{i}\Bigr)+\operatorname{ReLU}\left(x_{i}-w_{i}^{+}\right)  +\operatorname{ReLU}\Bigl(\frac{h_{i}}{2}-y_{i}\Bigr)+\operatorname{ReLU}\left(y_{i}-h_{i}^{+}\right)\big].
\]
Since $f_0$ and $f_b$ are sums of finite numbers of piecewise-linear functions, their subgradients are bounded; that is,
 there exists $G>0,B>0$ such that
 \begin{equation}
 \|v\|\le L,\ \|s\|\le B, \quad \text{for all } v\in\partial f_0(z), s\in\partial f_b(z).
 \end{equation}
Define $\mu := \min_{(x_i,y_i)\in V} \gamma_i$. To make the analysis easier, we assume that
\begin{equation}\label{eq:penalty-condition}
\mu > L.
\end{equation}

The SGD update can be further expressed as:
\[
z_{k+1}=z_k-\alpha_k \hat{g}_k,\qquad
\hat{g}_k=v_k+s_k+\xi_k,\quad v_k\in\partial f_0(z_k),\ s_k\in\partial f_b(z_k),
\]
with Gaussian noise $\xi_k\sim\mathcal N(0,\Sigma_k)$ satisfying
\[
\mathbb{E}[\xi_k\mid \mathcal F_k]=0,\qquad
\mathbb{E}[\|\xi_k\|^2\mid \mathcal F_k]\le \sigma^2 \ \ \text{a.s.}.
\]
And we assume the stepsize satisfies:
\begin{equation}\label{eq:stepsize-condition}
\alpha_k>0, \quad \lim_{k \to \infty}\alpha_k=0, \quad \sum_{k=0}^\infty \alpha_k=\infty, \quad \sum_{k=0}^\infty \alpha_k^2<\infty.
\end{equation}

To prove the boundedness of the sequence generated by the algorithm, we need to use the following Robbins–Siegmund almost-supermartingale theorem.

\begin{lemma}[Robbins-Siegmund~\cite{robbins1971convergence}]\label{lemma1}
Let $(X_k)$, $(A_k)$, $(B_k)$ and $(C_k)$ be sequences of nonnegative integrable random variables
 on some arbitrary probability space and adapted to the filtration $(\mathcal{F}_k)$, 
 with $\sum_{i=0}^\infty A_i < \infty$ and $\sum_{i=0}^\infty C_i < \infty$ almost surely and
\[
\mathbb{E}[X_{k+1} \mid \mathcal{F}_k] \leq (1 + A_k) X_k - B_k + C_k
\]
almost surely for all $k \in \mathbb{N}$. Then, almost surely, $(X_k)$ converges and $\sum_{i=0}^\infty B_i < \infty$.
\end{lemma}

\begin{theorem}\label{zkbound}
%Let $z=(x_1,y_1,\ldots,x_n,y_n)\in\mathbb{R}^{2n}$, 
For each block $i$, define the admissible ranges
\[
L_{x,i}=\tfrac{w_i}{2},\quad U_{x,i}=W-\tfrac{w_i}{2},\qquad
L_{y,i}=\tfrac{h_i}{2},\quad U_{y,i}=H-\tfrac{h_i}{2}.
\]
Set the n-dimensional box
\[
Q=\prod_{i=1}^n [L_{x,i},U_{x,i}]\times [L_{y,i},U_{y,i}].\]
Let $p(z)$ denote the Euclidean projection onto the box $Q$, and $d(z)$ represent the projection distance, i.e.
\[p(z) :=\arg\min_{y\in Q}\|z-y\|_2, \quad
d(z)=\mathrm{dist}(z,Q):=\|z-p(z)\|_2.
\]
Let $(\mathcal{F}_k)$ be the natural filtration generated by 
$ \{\eta_0, g_0, \eta_1, g_1, \ldots, \eta_k, g_k \} $,
i.e., all randomness observed up to step~$k$. Assume the mini-batch estimator is conditionally unbiased $(i.e.\,\mathbb{E}[g_k|\mathcal F_k]\in\partial F(z_k))$.
Then
\begin{enumerate}
\item[(i)]
For every $z$ and $s\in\partial f_b(z)$, we have
\begin{equation}\label{eq:inward}
\langle z-p(z),\, s\rangle \ \ge\ \mu\, d(z).
\end{equation}

\item[(ii)]
Let $d_k=d(z_k)$ and $p_k=p(z_k)$, where $\{z_k\}$ are the iterative variables generated by Algorithm~\ref{RBSM} with conditions \eqref{eq:penalty-condition}, \eqref{eq:stepsize-condition}. Then
\begin{equation}\label{eq:recursion-noisy}
\mathbb{E}\!\left[d_{k+1}^2\mid \mathcal F_k\right]
\ \le\ d_k^2 \;-\; 2(\mu-L)\,\alpha_k\, d_k \;+\; \alpha_k^2\, C^2,
\end{equation}
for some finite constant $C>0$.

\item[(iii)]
If $\mu\ge L$, then $\{z_k\}$ is almost surely bounded; more precisely, $d(z_k)\to 0$ almost surely.
\end{enumerate}
\end{theorem}

\begin{proof}
(i) It suffices to prove the claim for each coordinate and then sum.
Fix $i$ and consider the $x_i$-coordinate; the $y_i$-case is identical.
Let
$
p_{x,i}=\min\{\max\{x_i,L_{x,i}\},\,U_{x,i}\}
$
be the clipping of $x_i$ to the interval $[L_{x,i},U_{x,i}]$.
From the definition of the one-dimensional penalty
$\ell_{b,x,i}(x_i)=\gamma_i\big(\ReLU(L_{x,i}-x_i) + \ReLU(x_i-U_{x,i})\big)$, we have
\[
s_{x,i}\in \partial \ell_{b,x,i}(x_i)=
\begin{cases}
\{-\gamma_i\},      & x_i<L_{x,i},\\[2pt]
[-\gamma_i,\,0],    & x_i=L_{x,i},\\[2pt]
\{0\},         & L_{x,i}<x_i<U_{x,i},\\[2pt]
[0,\,\gamma_i],     & x_i=U_{x,i},\\[2pt]
\{\gamma_i\},       & x_i>U_{x,i}.
\end{cases}
\]
We distinguish three cases.
\begin{enumerate}
\item[1)] 
$x_i>U_{x,i}$. Then $p_{x,i}=U_{x,i}$, $x_i-p_{x,i}>0$ and $s_{x,i}=\gamma_i$.
Hence $(x_i-p_{x,i})\,s_{x,i}=(x_i-U_{x,i})\gamma_i=\gamma_i\,|x_i-p_{x,i}|$.
\item[2)] 
$x_i<L_{x,i}$. Then $p_{x,i}=L_{x,i}$, $x_i-p_{x,i}<0$ and $s_{x,i}=-\gamma_i$.
Hence $(x_i-p_{x,i})\,s_{x,i}=(x_i-L_{x,i})(-\gamma_i)=\gamma_i\,|x_i-p_{x,i}|$.
\item[3)]
$L_{x,i}\le x_i\le U_{x,i}$. Then $p_{x,i}=x_i$ and $(x_i-p_{x,i})\,s_{x,i}=0
= \gamma_i\,|x_i-p_{x,i}|$. 
\end{enumerate}
Therefore, for every coordinate,
\[
(x_i-p_{x,i})\,s_{x,i}\ \ge\ \mu\,|x_i-p_{x,i}|,\qquad
(y_i-p_{y,i})\,s_{y,i}\ \ge\ \mu\,|y_i-p_{y,i}|.
\]
Summing over $i=1,\ldots,n$ gives
\[
\langle z-p(z),\,s\rangle
=\sum_{i=1}^n\big[(x_i-p_{x,i})s_{x,i}+(y_i-p_{y,i})s_{y,i}\big]
\ \ge\ \mu\sum_{i=1}^n\big(|x_i-p_{x,i}|+|y_i-p_{y,i}|\big)
= \mu\,\|z-p(z)\|_1.
\]
Finally, $\|u\|_1\ge \|u\|_2$ for all $u$, hence
\[
\langle z-p(z),\,s\rangle\ \ge\ \mu\,\|z-p(z)\|_2
=\mu\, d(z),
\]
which is exactly \eqref{eq:inward}.

\smallskip
{(ii)}
From the properties of projection, we know
\[
d(z_{k+1}) = \|z_{k+1} - p(z_{k+1})\| \leq \|z_{k+1} - p(z_k)\| = \|z_k - \alpha_k \hat{g}_k - p(z_k)\|.
\]
Squaring and expanding the above inequality to get
\[
d_{k+1}^2 \le d_k^2 -2\alpha_k\langle z_k-p_k, \hat{g}_k\rangle + \alpha_k^2\|\hat{g}_k\|^2.
\]
Take conditional expectation and decompose $\hat{g}_k=v_k+s_k+\xi_k$:
\[
\mathbb{E}[d_{k+1}^2\mid \mathcal F_k]
\le d_k^2 -2\alpha_k\langle z_k-p_k, v_k+s_k\rangle
+ \alpha_k^2\,\mathbb{E}[\|v_k+s_k+\xi_k\|^2\mid \mathcal F_k].
\]
By Cauchy-Schwarz inequality and equation \eqref{eq:inward}, we have
\[
\begin{aligned}
\mathbb{E}[d_{k+1}^2\mid \mathcal F_k]
&\le d_k^2 -2\alpha_k\langle z_k-p_k, v_k\rangle-2\alpha_k\langle z_k-p_k, s_k\rangle
+ \alpha_k^2\,\mathbb{E}[\|v_k+s_k+\xi_k\|^2\mid \mathcal F_k],\\
&\le d_k^2 +2\alpha_k d_k\|v_k\|-2\alpha_k \mu d_k
+ \alpha_k^2\,\mathbb{E}[\|v_k+s_k+\xi_k\|^2\mid \mathcal F_k],\\
&\le d_k^2 -2(\mu-L)\alpha_k  d_k
+ \alpha_k^2\,\mathbb{E}[\|v_k+s_k+\xi_k\|^2\mid \mathcal F_k].
\end{aligned}
\]
Using $(a+b+c)^2\le 3(\|a\|^2+\|b\|^2+\|c\|^2)$, $\|v_k\|\le L$,
$\|s_k\|\le B$, and
$\mathbb{E}[\|\xi_k\|^2\mid\mathcal F_k]\le \sigma^2$ gives
$\mathbb{E}[\|\hat{g}_k\|^2\mid\mathcal F_k]\le C^2$, where $C^2=3(L^2+B^2+\sigma^2)$. Hence 
\[
\mathbb{E}\!\left[d_{k+1}^2\mid \mathcal F_k\right]
\ \le\ d_k^2 \;-\; 2(\mu-L)\,\alpha_k\, d_k \;+\; \alpha_k^2\, C^2.
\]

\textbf{(iii) }
Let $X_k=d_k^2$, $B_k=2(\mu-L)\alpha_k d_k$, $C_k=\alpha_k^2 C^2$.
Then
\[\mathbb{E}[X_{k+1}\mid \mathcal F_k]\le X_k-B_k+C_k.\]
Since $\sum_k\alpha_k^2<\infty$, $\mu> L$, 
the Robbins--Siegmund almost supermartingale theorem (Lemma~\ref{lemma1}) implies:
(1) $X_k$ (i.e. $d^2_k$) converges a.s.\ to a finite value;
(2) $\sum_k B_k<\infty$ almost surely.

From $\sum_k\alpha_k=\infty$, we know that
if $d_k \to a > 0$, then for sufficiently large $k$, $d_k \ge a/2$ almost surely.
Then $\sum_k B_k \ge \sum_k a(\mu-L)\alpha_k = \infty$ almost surely, which is a contradiction.
Thus, we have $d_k\to 0$ almost surely.
\end{proof}

Next, we prove that our algorithm converges to the Clarke-stationary set.

\begin{theorem}[Convergence to the Clarke-stationary set]\label{thm:conv-to-S}
Let $\{z_k\}$ be the optimization sequence generated by Algorithm~\ref{RBSM} with conditions \eqref{eq:penalty-condition}, \eqref{eq:stepsize-condition}.
Define $S=\{z:\ 0\in\partial f(z)\}$, where $\partial f(z)$ denotes the Clarke subdifferential.
Then $\mathrm{dist}(z_k,S)\to 0$ almost surely.
%\begin{itemize}
%\item[(i)] Every cluster point $z^\star$ of $(z_k)$ satisfies $0\in\partial^C F(z^\star)$.
%\item[(ii)]
%Consequently, $\mathrm{dist}(z_k,S)\to 0$.
%\end{itemize}
\end{theorem}

\begin{proof}
We first prove the deterministic case, i.e. $\epsilon_k=0$ in \eqref{eq:rbsm-one-step}. 
We assume that all affine regions are simplices; thus, they are convex. Otherwise, we can divide a nonsimplex region into several simplex regions by adding hyper-edges.
We split the remaining proof into several steps.

{(i)} %By Theorem~\ref{zkbound}, the sequence $(z_k)$ is a.s.\ bounded (for $\mu\ge L$).
By the definition of our objective, $S$ is bounded and closed. By Theorem \ref{gradnorm-bound}, $G_{S^c}$ has an upper bound $L_g$ and a lower bound $\beta$.
Moreover, for any $z\notin S$, there exists a unit vector $u_z$ and $\beta>0$ such that
\begin{equation}
\langle v,u_z\rangle\ \ge\ \beta \qquad  \forall\, v\in\partial f(z).\label{sepa}
\end{equation}

(ii) Define the open halfspace $H_\beta^u:=\{w:\langle w,u\rangle>\beta\}$.
By the closeness of $S$, We can show that for any $z\notin S$, there is a neighborhood $U_z$ of $z$, a unit vector $u_z$, such that 
$\partial f(z')\subset H_\beta^{u_z}$ for all $z'\in U_z$, which means that while $z_k$ is inside $U_z$, each iteration step will have a movement along $-u_z$ direction that is no less than $\beta$. Then, by the condition $\sum \alpha_k \to \infty$, the optimization series will exit in finite steps from any given neighborhood $U_z$ of finite size. That means
the deterministic subgradient descent method can't converge to a point $z^* \notin S$.

(iii) By Theorem~\ref{zkbound}, the sequence generated by Algorithm~\ref{RBSM} is a.s.\ bounded (for $\mu\ge L$). Similarly, the deterministic optimization sequence is bounded. Therefore, there exists at least one cluster point. If a cluster point $z^* \notin S$, then by (ii), we know that there exists a neighborhood $U(z^*;\delta)$ with 
$\delta < dist(z^*, S)$, and the sequence $\{z_k\}_{k=1}^\infty$ must enter $U(z^*;\delta)$ infinitely many times. Next, we want to show that when
$k$ is large enough, i.e., $\max_{k'>k}\alpha_k'$ is small enough, re-entering the neighborhood $U(z^*; \delta)$ is not possible. 

Since the sequence is bounded, we can cover the admissible region using a finite number of small balls
with a finite lower bound on the radius. We can take the radius of the balls small enough such that all the affine pieces in a small ball have a non-empty common intersection. We call the intersection ``central structure'', denoted by $C_s$, of the small ball $U_s$.
For simplicity, we can merge the balls that involve the same central structure into one set, which is still denoted by $U_s$. Next, we prove that starting from a given point $z_k$ near the center of $U_s$, for 
large enough $k$,  $f(z_k)-f(z_{k+M+1})\ge \delta >0$, where $z_{k+M+1}$ is the first point in the optimization sequence $\{ z_k, z_{k+1}, \ldots\}$ that is not inside $U_s$.
We discuss this by considering different types of central structures. Denote by $m(U_s)$ the number of affine pieces inside $U_s$.

a) If $m(U_s)=1$, then the central structure is a $2n$-dimensional (full dimension) element; i.e., this region $U_s$ is totally contained within one affine region. 
Then, as long as the optimization point is inside $U_s$, the objective function decays, except for the last point that exits $U_s$. 
Let $\delta_s :=dist(z_k, \partial {U_s})$.
Take $k$ large enough, such that $\alpha_{k'} < \beta \delta_s /L_g$ for all $k'>k$. Since $\sum \alpha_k \to \infty$, there exists a large enough $M$ such that
$\sum_{j=0}^{M-1} \alpha_{k+j} \beta >\delta_s$.
Then we have
\begin{equation}
f(z_k) - f(z_{k+M+1}) \ge \sum_{j=0}^{M-1} \alpha_{k+j} |g_s|^2 - \alpha_{k+M} L_g 
\ge \sum_{j=0}^{M-1} \alpha_{k+j} \beta^2 - \alpha_{k+M} L_g 
> 0.
\end{equation}

b) If $m(U_s)=2$, then the central structure $C_s$ of $U_s$ is a $2n-1$ dimensional interface (face).
Without loss of generality, we assume that $C_s$ is at the center of $U_s$.
A starting point $z_k$ close to the center of $U_s$ is also close to $C_s$.
For $n\ge 1$, $C_s$ is not a vertex. We can take the unit vector $u$ in the separating hyperplane argument parallel to $C_s$. This can be done because all the gradients projected into the subspace defined by $C_s$ are identical, by the continuity of $f$. In fact, $u$ can be taken as any normalized vector that is positively correlated to the gradient of $f$ confined on $C_s$. Denote by $g_i^u := \sum_{\ell=1}^{2n-1}\langle g_i, u_\ell\rangle u_\ell$, $g_i^\perp := g_i - g_i^u$, where $\{u_\ell, \ell=1,\ldots, 2n-1 \}$ is a set of orthonormal basis of $C_s$ with each axis positively correlated to the gradient of $f$ confined on $C_s$. We consider the optimization sequence projected on the hyper-plane that is orthogonal to $C_s$.
If $\langle g_1^\perp, g_2^\perp \rangle \ge 0 $, then the optimization sequence $\{z_{k+j}\}$ will not pass through $C_s$ more than once before exiting $U_s$. 
Since we assume the initial point $z_k$ is very close to $C_s$, we can expect that the pass-through occurs at the first optimization step, if it exists. 
Correspondingly, if $k$ is large enough (similar to the $m(U_s)=1$ case), the objective function will have a sufficient decay before exiting $U_s$ as
$$
f(z_k) - f(z_{k+M+1}) \ge \sum_{j=1}^{M-1} \alpha_{k+j} \beta^2 - \alpha_{k+M} L_g > 0.
$$
If $\langle g_1^\perp, g_2^\perp \rangle < 0 $, then for sufficiently large $k$, the optimization sequence $\{ z_{k+j} \}$ will either not pass through $C_s$
or pass through $C_s$ many times before exiting $U_s$. That is the reason why we call $C_s$ a central structure. 
In the first case, the objective function decreases at every time step, then the situation is similar to the case $m(U_s)=1$. In the second case, each time the sequence passes through $C_s$, the objective function will experience a jump but then decrease towards $C_s$.  
Consider the objective values at $z_k, z_k^*, z_{k+M}^*, z_{k+M}, z_{k+M+1}$, 
where $z_k^*, z_{k+M}^*$ denotes the first and last intersection points of the optimization sequence with $C_s$, respectively. 
We have
$$
f(z_k) - f(z_k^*) \ge 0,
\quad 
f(z_k^*) - f(z_{k+M}^*) \ge \sum_{j=1}^{M-1} \alpha_{k+j} \beta,
$$
and 
$$
 f(z_{k+M}^*) - f(z_{k+M}) \ge - \alpha_{k+M'} L_g,
\quad
 f(z_{k+M}) - f(z_{k+M+1}) \ge - \alpha_{k+M} L_g,
$$
where $M'$ is a number between 0 and $M$.
So, by taking $k$ large enough, we have
$$
f(z_k) - f(z_{k+M+1}) \ge \sum_{j=1}^{M-1} \alpha_{k+j} \beta - (\alpha_{k+M}+\alpha_{k+M'}) L_g > 0.
$$

c) Now, we consider the case $m(U_s)>2$ and $\text{dim}(C_s)=0$, i.e., $C_s$ is a vertex. 
Since $U_s\cap S=\emptyset$, by our choice of $U_s$ and the separating hyperplane argument, there exists a unit vector $u_s$ such that
\begin{equation}
\langle v, u_s\rangle \ge \beta > 0, \ \forall v \in \partial f(z),\, z\in U_s.
\end{equation}
That means when $z_k$ is close to $C_s$, the optimization point will move a certain distance along the $u_s$ direction. In other words, the optimization sequence $\{ z_{k+j} \}$ will pass through the hyperplane that is centered at $C_s$ and orthogonal to $u_s$ at most once. After passing-through, the sequence $\{ z_{k+j} \}$ will remain inside one affine region or linger around another central structure $C_{s'}$ with $\text{dim}>0$ and linked to $C_s$. For the case $n=1$, this leads to the situations already considered in a) and b), by which arguments we can take large enough $k$ to make the objective function have sufficient decay during the entire optimization procedure inside $U_s$ and $U_{s'}$.
For the case $n>1$, we provide an analysis of lingering around another central structure $C_{s'}$, with $\text{dim}>0$ next.

d) Now, we consider the case $m(U_s)>2$ and $d:=\text{dim}(C_s)>0$. In this case, as done in b), the $u$ vector in the separating hyperplane argument can be taken as any normalized vector that is positively correlated with the gradient of $f$ confined to $C_s$. Denote by $g_i^u := \sum_{\ell=1}^{d}\langle g_i, u_\ell\rangle u_\ell$, $g_i^\perp := g_i - g_i^u$, where $\{u_\ell, \ell=1,\ldots, d \}$ is a set of orthonormal bases of $C_s$, each axis is positively correlated with the gradient of $f$ confined to $C_s$. 
We consider the optimization sequence projected onto the hyperplane that is orthogonal to $C_s$.
If there exists a unit vector $u_s^\perp \perp u_s$ and $\beta_s^\perp > 0$ such that
 $$\langle v, u_s^\perp \rangle \ge \beta_s^\perp,\quad \forall v\in \clconv\{g_i^\perp, i=1,\ldots, d \}, $$ 
 then the optimization sequence $\{ z_{k+j}\}$ will not pass through the hyperplane $ (u_s^\perp)^t z = c$ (which contains $C_s$) more than once, which means the optimization sequence will move into another neighborhood $U_{s'}$ containing a central structure $C_{s'}$ with a co-dimension smaller than $C_s$'s. 
 Now, consider the case where the optimization sequence lingers about $C_s$ for many iterations, for which there must be $0 \in \clconv\{g_i^\perp, i=1,\ldots, d \}$.
 We only need to consider the case that $dist(z_{k+j},C_s)< c \max_j{\tau_{k+j}} L_g /\beta =: \delta_k$, 
 where $c$ is a constant related to the minimum angle among the faces that are connected to $C_s$, because in this case, after a jump in the objective function, there will be a large enough affine region to let the objective function decay the jumped value.
 Now, consider the objective values at $z_k, z_k^u, z_{k+M}^u, z_{k+M}, z_{k+M+1}$, where $z_k^u$ and $z_{k+M}^u$ are the projections of $z_k$ and $z_{k+M}$ onto $C_s$.
 We have
$$
f(z_k^u) - f(z_{k+M}^u) \ge \sum_{j=0}^{M-1} \alpha_{k+j} \beta^2,
\quad |f(z_k) - f(z_k^u)| \le \delta_k L_g,
$$
and 
$$
| f(z_{k+M}^u) - f(z_{k+M})| \le \delta_k L_g,
\quad
| f(z_{k+M}) - f(z_{k+M+1})| \le \alpha_{k+M} L_g.
$$
So, by taking $k$ large enough, we have
$$
f(z_k) - f(z_{k+M+1}) \ge \sum_{j=0}^{M-1} \alpha_{k+j} \beta - (2\delta_k+\alpha_{k+M}) L_g > 0.
$$

Since we can cover the entire admissible domain by a finite number of central structure neighborhoods $U_s, s=1,\ldots, n_s$. By the above argument, starting from a point $z_k$ near the center of any given $U_s$, with $k$ sufficiently large, the objective value at the first point of exiting $U_s$ decreases.
Then, for a cluster point $z^* \notin S$, once $z_k$ exits $U(z^*; \delta)$, the optimization sequence $\{ z_{k+j} \}$ with $k$ sufficiently large will not re-enter the neighborhood of $z^*$ since any path returning to the neighborhood of $U(z^*;\delta'), \delta' < \delta$ will pass through a certain number of central structure neighborhoods $U_s, s=1,\ldots, n_{s,z}$. After that, the objective function value will be smaller than the minimum function value on $U(z^*;\delta')$. Thus, it can't re-enter $U(z^*; \delta')$.
We proved that $dist(z_k, S)\to 0$, as $k\to \infty$ for the deterministic case.

For the stochastic case, by taking the expectation of both sides of \eqref{eq:rbsm-one-step}, we obtain
\begin{equation}
\mathbb{E}(z_{k+1} \mid \mathcal{F}_k) = z_k - \alpha_k  g(z_k).
\end{equation}
At each step, the mean of $z_k$ follows the same dynamics as in the deterministic case.  Thus,
$dist(\mathbb{E}(z_k), S) \to 0$, as $k\to \infty$.
Next, we estimate the variance of the iteration
\begin{equation}
\begin{aligned}
\text{Var}(z_{k+1}) &= \text{Var}(z_k - \alpha_k g_k) + \alpha_k^2 \text{Var}(\xi_k)\\
    & = \text{Var}(z_k) + \alpha_k^2 \text{Var}(g_k) - 2\alpha_k \mathbb{E}\bigl((z_k - \mathbb{E} z_k)(g_k - \mathbb{E} g_k)\bigr) + \alpha_k^2 \text{Var}(\xi_k).
\end{aligned}
\end{equation}
When $k$ is large enough such that $\mathbb{E}z_k$ is very close to $S$, we have
$z_k-\mathbb{E} z_k$ and $g_k - \mathbb{E} g_k$ that are nonnegatively related, i.e.,
$\mathbb{E}\bigl((z_k - \mathbb{E} z_k)(g_k - \mathbb{E} g_k)\bigr)\ge 0$. Therefore
\begin{equation}
\text{Var}(z_{k+1}) \le 
    \text{Var}(z_k) + \alpha_k^2 L_g^2  + \alpha_k^2 \text{Var}(\xi_k),
\end{equation}
for sufficiently large $k$. Summing up the above equation for $k=j, j+1, \ldots, j+m$, we get
\begin{equation}
    \text{Var}(z_{j+m+1}) \le \text{Var}(z_j) + \sum_{k=j}^{j+m}\alpha_j^2 L_g^2  + \sum_{k=j}^{j+m} \alpha_j^2 \text{Var}(\xi_j).
\end{equation}
Note that when calculating the above inequality with the given filtration $\mathcal{F}_j$, we have $\text{Var}(z_j|\mathcal{F}_j ) = 0$.
Then, due to the fact that $\text{Var}(\xi_j)$ and $\sum_{j=1}^\infty \alpha_j^2$ are bounded,
we obtain that
$\lim_{k\to \infty} \max_{k'>k} \text{Var}(z_{k'}|\mathcal{F}_k) =0$, i.e. $\mathrm{dist}(z_k,S)\to 0$ almost surely. The theorem is proved.
\end{proof}

\section{ Numerical experiments}\label{sec4}
\subsection{Experimental setup}

In this section, we present a series of numerical experiments to assess the performance of three different methods: gradient descent (GD)~\cite{nocedal1999numerical}, Adam~\cite{kingma2014adam}, and RBSM. The evaluation is conducted using the GSRC benchmarks, which are widely used to test and evaluate electronic design automation tools and algorithms. The code is written in Python and run on a standard PC with an i9-13980HX microprocessor and 32GB of memory.

\subsection{Results on GSRC benchmark}
We implemented the RBSM method, GD method, and Adam method in Python. The GD method only uses the built-in SGD subroutine of Torch for optimization, while the Adam method is combined with a series of proposed enhanced techniques, including parameter updates, mean-field force, and weighted sampling. This is to compare the optimization algorithm and enhanced techniques from two perspectives. In the RBSM and Adam methods, the parameters are set to identical values: average field force coefficient $\alpha=5$, batch size is $1/5$ of the total number of nets, and the initial value of parameter $\gamma$ is $1000$. In the GD method, $\gamma$ is always set to a constant of $1000$. We verify the feasibility and performance of the proposed algorithm by GSRC plane planning benchmarks. The final GSRC graphic design results were used as a test. We compared the wirelength and overlap of the different methods, as well as the visual effects of the placement(see Appendix). Table \ref{tab:1} shows the experimental data of GSRC and the results of different algorithms, while Table \ref{tab:compare1} presents the placement and legalization results of different algorithms in the benchmark test.

\begin{table}[ht]
\centering
\caption{GSRC benchmark circuit standards}
\label{tab:1}
    \begin{tabular}{lrrrrrr}
        \toprule
        circuit & module & terminal & net  & pin  & \textbf{W} $\times$ \textbf{H} & Die-Area$(\mathrm{mm}^2)$  \\
        \midrule
        n10     & 10     & 69       & 118  & 248  & $800\times800$                 & 640000                     \\
        n30     & 30     & 212      & 349  & 743  & $800\times800$                 & 640000                     \\
        n50     & 50     & 209      & 485  & 1050 & $800\times800$                 & 640000                     \\
        n100    & 100    & 334      & 885  & 1873 & $800\times800$                 & 640000                     \\
        n200    & 200    & 564      & 1585 & 3599 & $800\times800$                 & 640000                     \\
        n300    & 300    & 569      & 1893 & 4358 & $800\times800$               & 640000                     \\
        \bottomrule
    \end{tabular}
\end{table}

\begin{table}[ht]
    \centering
    \caption{The experimental results of RBSM, GD, and Adam methods in different circuits of the GSRC benchmark. Where LHPWL means the wirelength after legalization, percentage means the overlap percentage relative to the total cell area. In the RBSM and Adam methods, the average field force coefficient $\alpha=5$, the batch size is $1/5$ of the total number of nets, and the initial value of parameter $\gamma$ is $1000$. In the GD method, $\gamma$ is set to a constant of $1000$. }
    \label{tab:compare1}
    {
         \begin{tabular}{l|rr|rrrr}
                \toprule
                & \multicolumn{2}{c|}{GSRC} & \multicolumn{4}{c}{RBSM} \\
                \hline
                Circuits & HPWL & overlap(percentage) & HPWL & overlap & time(s) & LHPWL \\
                \midrule
                n10 & 64299 & 0 & 56902 & 7206(3.25\%) & 7.48 & 56894 \\
                n30 & 179811 & 0 & 156761 & 1896(0.91\%) & 6.89 & 157261 \\
                n50 & 234281 & 0 & 199356 & 602(0.30\%)  & 3.30 & 199236 \\
                n100 & 395719 & 0 & 328705 & 1470(0.84\%)  & 5.78 & 328991 \\
                n200 & 738707 & 0 & 571720 & 2450(1.39\%) & 17.76 & 574778 \\
                n300 & 937608 & 0 & 694527 & 3634(1.33\%) & 36.12 & 698867.6  \\
                \bottomrule
            \end{tabular}
            
         \begin{tabular}{l|rrrr|rrrr}
               \toprule
            & \multicolumn{4}{c|}{GD} & \multicolumn{4}{c}{Adam} \\
            \midrule
            Circuits & HPWL & overlap(\%) & time(s) & LHPWL & HPWL & overlap(\%) & time(s) & LHPWL \\
            \midrule
            n10 & 63959 & 82(0.037\%) & 0.45 & 63955 & 59631 & 29460(13.29\%) & 8.83 & 62240 \\
            n30 & 188770 & 612(0.29\%) & 0.40 & 188192 & 173825 & 35613(17.07\%) & 9.64 & 185765 \\
            n50 & 224062 & 1(0.0005\%) & 0.63 & 224038 & 217210 & 42382(21.35\%) & 5.69 & 230797 \\
            n100 & 366449 & 1(0.0006\%) & 0.92 & 366344 & 376827 & 20581(11.70\%) & 8.31 & 389947 \\
            n200 & 728562 & 0(0\%) & 2.16 & 728475 & 697857 & 19104(10.87\%) & 20.18 & 720499 \\
            n300 & 978937 & 60(0.022\%) & 3.76 & 978438 & 797910 & 71746(26.26\%) & 47.87 & 914895 \\
            \bottomrule
            \end{tabular}}
\end{table}

Through the experimental results, we can observe the following points:
\begin{enumerate}
\item
The proposed placement algorithm effectively reduces overlap while significantly minimizing wirelength. This dual improvement demonstrates the method's strong potential, as conventional placement algorithms often struggle to simultaneously achieve these objectives during the global placement phase. Consequently, our approach facilitates more efficient legalization and detailed placement in subsequent stages.
\item
Compared to the GD algorithm without any enhanced techniques, the RBSM algorithm achieves a substantial reduction in wirelength. Notably, the wirelength after the legalization phase remains significantly shorter than that of the GD-based approach, underscoring the critical role of the enhancement techniques integrated into the RBSM algorithm.
\item
When the SGD optimizer in the RBSM algorithm is replaced with the Adam optimizer, the algorithm's performance deteriorates significantly, even falling below the performance of the GD method without enhancements. This indicates that gradient descent is more suitable for solving the rectified linear penalty model of placement problems due to its stability and convergence properties in this context.
\item
The proposed legalization algorithm not only eliminates overlap but also preserves the wirelength advantages achieved during placement. In some cases, it further reduces the wirelength. This result highlights the ability of the RBSM algorithm to maintain and even enhance placement quality, demonstrating its superiority in achieving short wirelength while ensuring overlap removal. 
\end{enumerate}

\subsection{Analysis of our enhanced technique}
In order to further illustrate the processing effectiveness of our optimization algorithm, we conducted a series of comparative experiments to assess the impact of four key enhancement techniques: degree-weighted sampling, adaptive parameter updating, mean-field force, and random disturbance on the algorithm, respectively. The results presented are the average outcomes from five independent runs of each algorithm. For clarity, each comparison isolates the effect of a single enhancement by modifying only the corresponding component of the algorithm, allowing for a precise evaluation of its contribution.

\begin{table}[ht]
    \caption{Comparative experiments on different enhancement techniques of RBSM method, with Original algorithm representing the results of RBSM method on corresponding circuits; Random batch represents the result of replacing weighted sampling with random sampling in the RBSM method (with the same sampling size); Fix gamma represents the result of using a fixed constant $\gamma=10000$ instead of adaptive $\gamma$; Without mean-field force means that mean-field force is not applied; Without random perturbation means the result without random perturbation.}
\label{tab:compare}
    \centering
  %\resizebox{\textwidth}{!}
 % {
\begin{tabular}{l|rrr|rrr|rrr}
\toprule
\multicolumn{1}{c|}{Circuit} & \multicolumn{3}{c|}{n100} & \multicolumn{3}{c|}{n200} & \multicolumn{3}{c}{n300} \\
  & HPWL & overlap & time & HPWL & overlap & time & HPWL & overlap & time   \\
\midrule
RBSM & 328705 & 1470 & 5.78 & 571720 & 2450 & 17.66 & 694527 & 3634 & 36.12 \\
 Random batch & 308666 & 3500 & 9.67 & 553329 & 4014 & 17.81 & 687753 & 5288 & 37.97 \\
 Fix gamma & 368020 & 1691 & 2.52 & 642329 & 1349 & 6.34 & 871881 & 2229 & 16.18 \\
No mean force & 332958 & 2099 & 5.75 & 622479 & 1629 & 14.99 & 767695 & 3317 & 31.99 \\
 No perturbation & 302856 & 3575 & 7.97 & 539778 & 7241 & 18.49 & 653182 & 10041 & 40.55 \\
\bottomrule
\end{tabular}
%}
\end{table}

1) {\textit{Degree-Weighted Sampling:}} We compared random sampling with degree-weighted sampling. The experimental results for the n100, n200, and n300 circuits, under the same number of iteration steps, demonstrate that the overlap for random sampling is one to two times higher than that for degree-weighted sampling. This finding highlights that degree-weighted sampling enhances convergence performance by allocating more computational resources to critical nets, thereby accelerating the optimization process.

2) {\textit{Adaptive Parameter Adjustment:}} The experimental results reveal that adaptive parameter updating significantly reduces wirelength compared to fixed-parameter algorithms. As the problem size increases, the performance gap becomes more pronounced. For instance, in the n100 circuit, the adaptive parameter adjustment reduces wirelength by $10\%$, and in the n300 circuit, it achieves a reduction of $20\%$. This demonstrates that adaptive parameter selection provides stronger adaptability and optimization capability for circuit placement, enabling significant reductions in wirelength and overlap. However, it is important to note that the runtime of the adaptive parameter algorithm is approximately twice that of the fixed-parameter algorithm, indicating potential areas for future optimization.

3) {\textit{mean-field force:}} As can be seen in the table, removing the mean-field force from the algorithm results in a $1\%$, $8\%$, and $10\%$ increase in wirelength for the n100, n200, and n300 circuits, respectively, compared to our full algorithm. This shows that adding the force of the aggregation effect can significantly improve the placement effect, and the effect becomes increasingly significant with the increase of the problem scale.

4) {\textit{Random Disturbance:}} Under a fixed number of iterations, incorporating random disturbances effectively reduces overlap. This is because random perturbations enable the algorithm to escape local minima, thereby accelerating convergence and improving overall optimization efficiency. This mechanism allows the algorithm to better explore the solution space, leading to improved placement quality.

\medskip
These experimental results collectively demonstrate the importance and effectiveness of the proposed enhancements in improving the performance of the placement algorithm. Each component contributes uniquely to addressing specific challenges in circuit placement.

\subsection{Results on ISPD2005 benchmark}
To further evaluate the scalability and effectiveness of our proposed method, 
we conduct experiments on large-scale ISPD2005 benchmark circuits and directly 
compare the results of our RBSM algorithm with DREAMPlace, a state-of-the-art 
GPU-accelerated global placement tool. 

Table \ref{tab:adaptec_stats} presents the ISPD2005 benchmark details, and Table \ref{tab:compare2} presents a detailed comparison between DREAMPlace 
and RBSM across four large benchmarks. Note that all results are reported before legalization, i.e., they correspond to the global placement stage only.
\begin{table}[ht]
\centering
\caption{Statistics of ISPD2005 benchmarks}
\label{tab:adaptec_stats}
\begin{tabular}{lrrrrrrr}
\toprule
Circuit & \makecell{Total\\Objects} & \makecell{Mov.\\Objects} & \makecell{Fixed\\Objects} 
        & Nets & \makecell{Total\\Pins} & \makecell{Pins.Mov.\\Objs} & \makecell{Pins.Fixed.\\Objs} \\
\midrule
adaptec1 & 211447 & 210904 & 543  & 221142 & 944053  & 923513  & 20540 \\
adaptec2 & 255023 & 254457 & 566  & 266009 & 1069482 & 1045699 & 23783 \\
adaptec3 & 451650 & 450927 & 723  & 466758 & 1875039 & 1843852 & 31187 \\
adaptec4 & 496045 & 494716 & 1329 & 515951 & 1912420 & 1876563 & 35857 \\
\bottomrule
\end{tabular}
\end{table}

\begin{table}[ht]
\centering
\caption{Comparison between DREAMPlace and our RBSM algorithm on ISPD2005 benchmarks. 
In this table, column \textit{overflow} refers to the density overflow percentage relative to the placement region capacity, 
whereas column \textit{overlap} represents the overlap percentage relative to the total cell area.}
\label{tab:compare2}
\begin{tabular}{l|rrr|rrrr}
\toprule
\multirow{2}{*}{Circuit} & \multicolumn{3}{c|}{DREAMPlace} & \multicolumn{4}{c}{RBSM} \\
 & HPWL  & overflow & epoch & HPWL & overflow & epoch & overlap \\
\midrule
adaptec1 & $7.03\times 10^7$ & $6.92\%$ & 607 & $5.05\times 10^7$ & $26.8\%$ & 85  & 10.36\% \\
adaptec2 & $7.93\times 10^7$ & $6.89\%$ & 637 & $8.06\times 10^7$ & $16.9\%$ & 100 & 11.67\% \\
adaptec3 & $1.86\times 10^8$ & $6.92\%$ & 679 & $1.92\times 10^8$ & $12.9\%$ & 320 & 4.74\% \\
adaptec4 & $1.69\times 10^8$ & $6.79\%$ & 697 & $1.79\times 10^8$ & $14.0\%$ & 220 & 3.25\% \\
\bottomrule
\end{tabular}
\end{table}

From Table \ref{tab:compare2}, we observe that the proposed RBSM algorithm 
achieves wirelength (HPWL) results comparable to DREAMPlace, while obtaining 
significantly smaller overlap (within $12\%$). This demonstrates that our method 
is effective in directly reducing both wirelength and overlap, which are the 
primary objectives in global placement. 

Although the density overflow and runtime of RBSM are relatively large, this outcome is expected because density optimization has not been explicitly incorporated into the algorithm design, and the current implementation has not been fully optimized for efficiency. Under these settings, such results are acceptable and reasonable. At the same time, RBSM already demonstrates strong competitiveness in the key placement objectives of wirelength reduction. Furthermore, since the number of epochs required by RBSM is actually smaller than that of DREAMPlace, we believe that with GPU parallelization and density optimization incorporated, RBSM has the potential to achieve runtime and overall placement quality comparable to or even better than DREAMPlace.

\section{Conclusion}\label{sec5}
In this paper, we take a new perspective to solve the placement problem and transform it into a neural net training problem. A nonsmooth optimization model for VLSI global placement is proposed, with half-perimeter wirelength minimization and cell overlap elimination as the primary objectives and constraints. A penalty function model is established as a loss function, which can be formulated as a deep neural network with only ReLU layers, without gradient vanishing and exploding issues. A stochastic subgradient method with operator splitting is proposed to solve the unconstrained nonsmooth optimization problem, and several optimization acceleration techniques are adopted to accelerate convergence. In addition, we provide a theoretical convergence analysis of the proposed algorithm, which guarantees its reliability for large-scale optimization. Numerical experiments on medium-scale GSRC benchmarks demonstrate that our model can effectively reduce HPWL and cell overlap, and the high-quality solutions obtained provide references for other imprecise algorithms. Moreover, large-scale experiments on ISPD2005 benchmarks show that our method achieves wirelength performance comparable to state-of-the-art tools with acceptable overlap. However, the current implementation is slower due to the lack of GPU acceleration and explicit density optimization. Improving efficiency and scalability will be the focus of future research.

\section*{Acknowledgements}
This work is partially supported by the National Key R\&D Program of China (grant no. 2021YFA1003601),
the National Natural Science Foundation of China (grant no. 92473208, 12494543, and 12171467), and
the Strategic Priority Research Program of the Chinese Academy of Sciences (grant no. XDA0480504). The computations were partially performed on the LSSC4 PC cluster of State Key Laboratory of Scientific and Engineering Computing.

\section*{Ethics declarations}

\subsection*{Conflict of interest}
The authors have no conflicts of interest to declare that are relevant to the content of this article.

\section*{Data availability statement}
The authors confirm that the data supporting the findings of this study are available within the article \cite{Gsrc} and its supplementary materials.

\begin{appendices}
    \section{ Visual placement effect display}
\begin{figure}[ht]
    \centering
    \includegraphics[width=\textwidth]{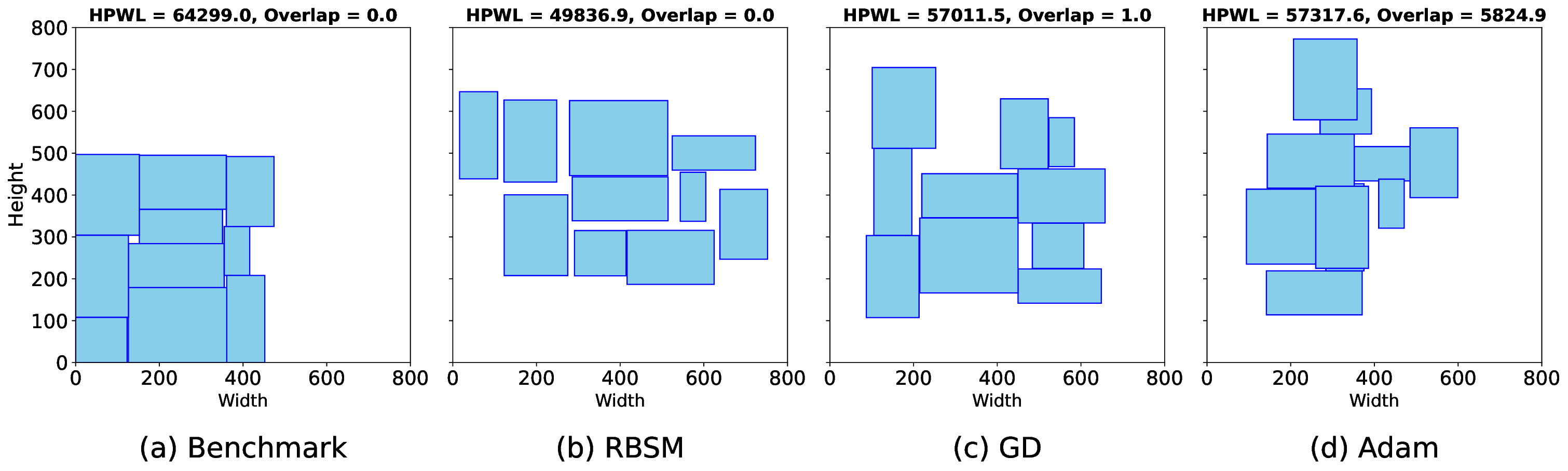}
    \caption{Placement comparison for $n = 10$ using different methods}
    \label{fig:n10_images2}
\end{figure}

\begin{figure}[ht]
    \centering
    \includegraphics[width=\textwidth]{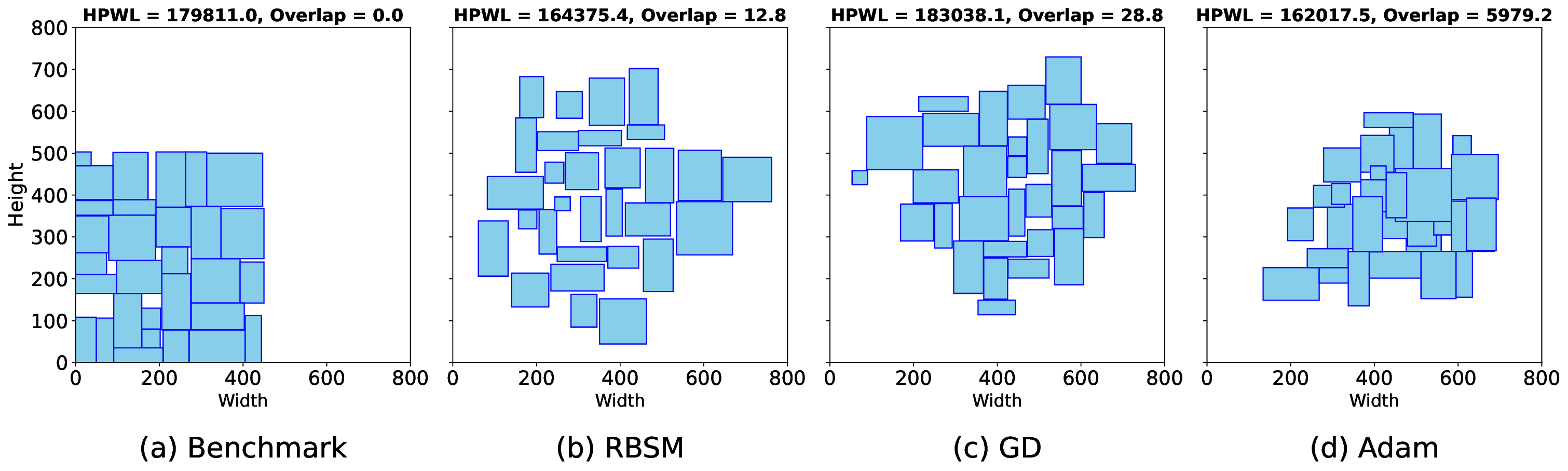}
    \caption{Placement comparison for $n = 30$ using different methods}
    \label{fig:n30_images2}
\end{figure}

\begin{figure}[ht]
    \centering
    \includegraphics[width=\textwidth]{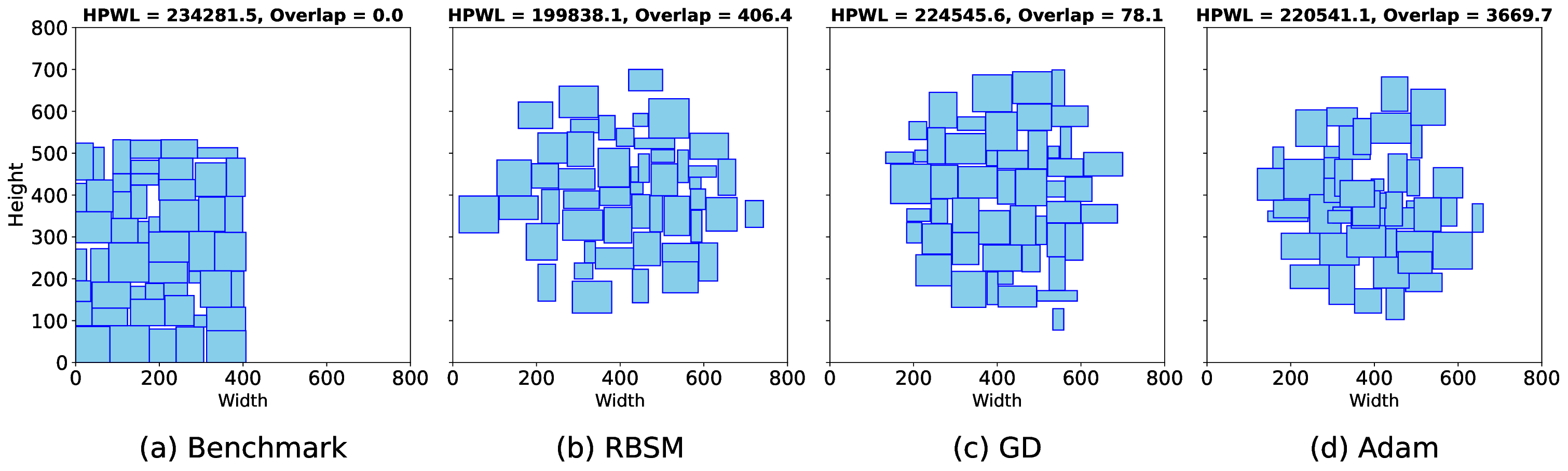}
    \caption{Placement comparison for $n = 50$ using different methods}
    \label{fig:n50_images2}
\end{figure}

\begin{figure}[ht]
    \centering
    \includegraphics[width=\textwidth]{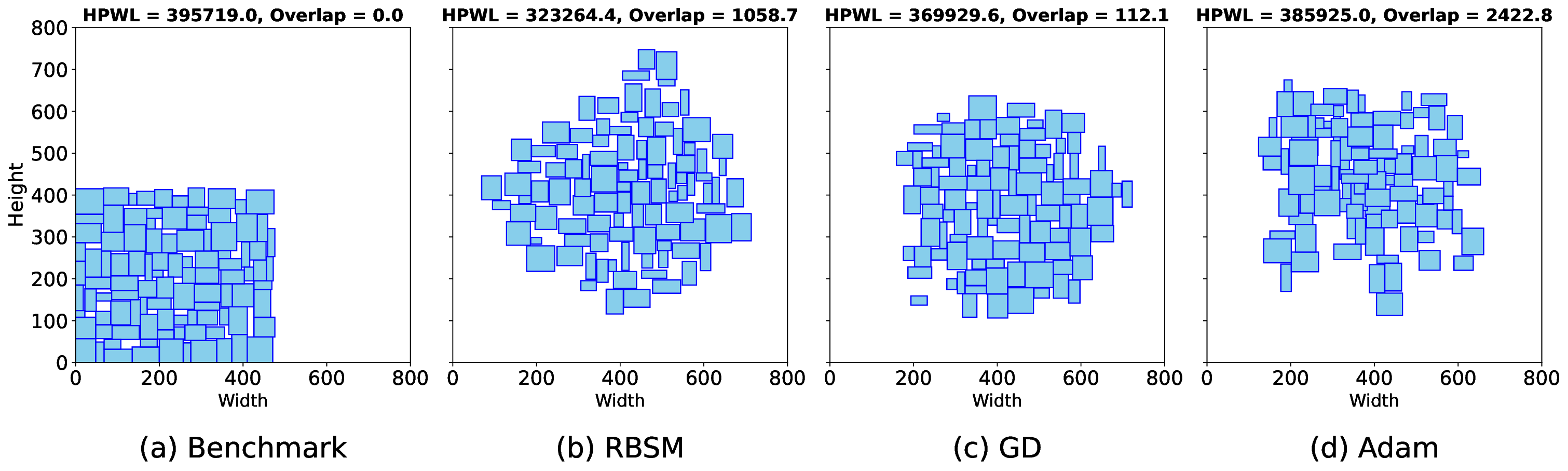}
    \caption{Placement comparison for $n = 100$ using different methods}
    \label{fig:n100_images2}
\end{figure}

\begin{figure}[ht]
    \centering
    \includegraphics[width=\textwidth]{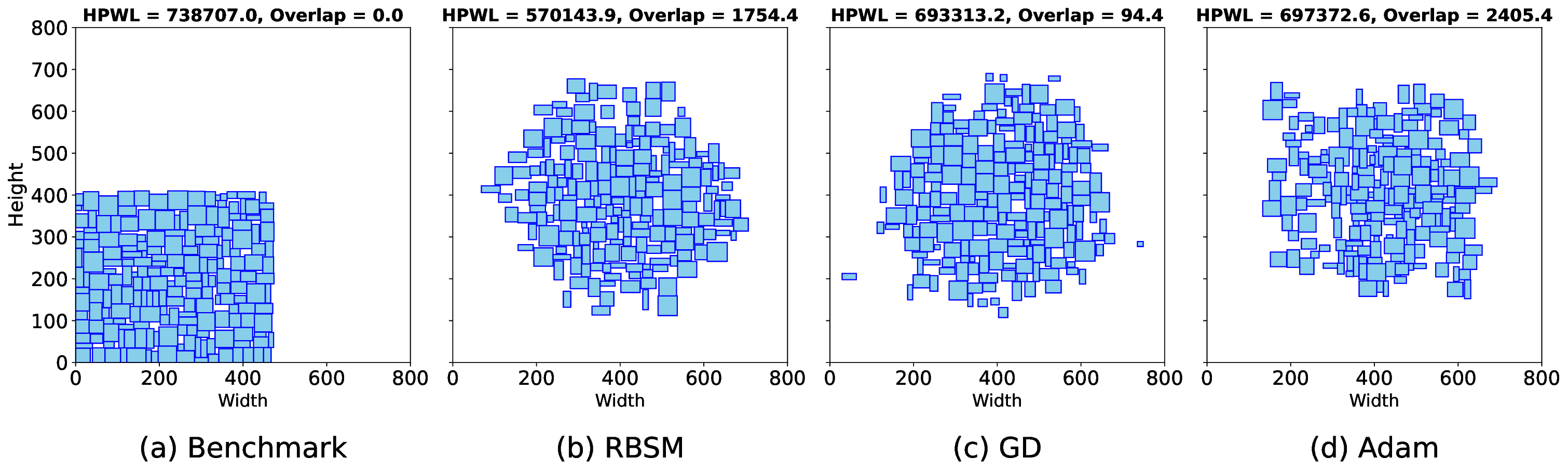}
    \caption{Placement comparison for $n = 200$ using different methods}
    \label{fig:n200_images2}
\end{figure}

\begin{figure}[ht]
    \centering
    \includegraphics[width=\textwidth]{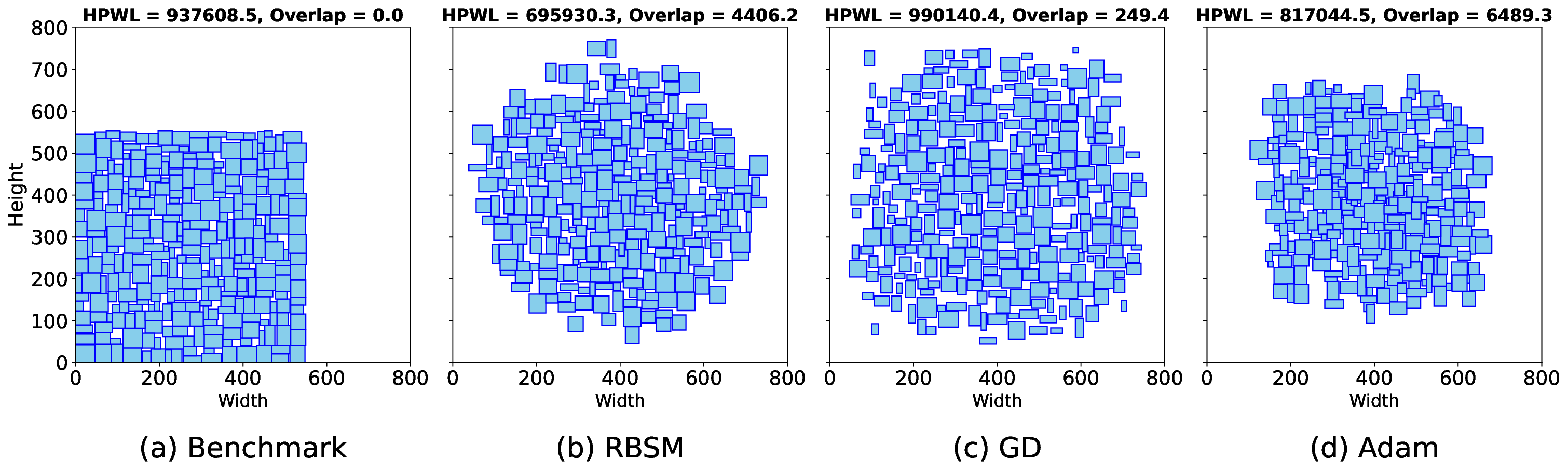}
    \caption{Placement comparison for $n = 300$ using different methods}
    \label{fig:n300_images2}
\end{figure}

\begin{figure}[ht]
    \centering
    \includegraphics[width=\textwidth]{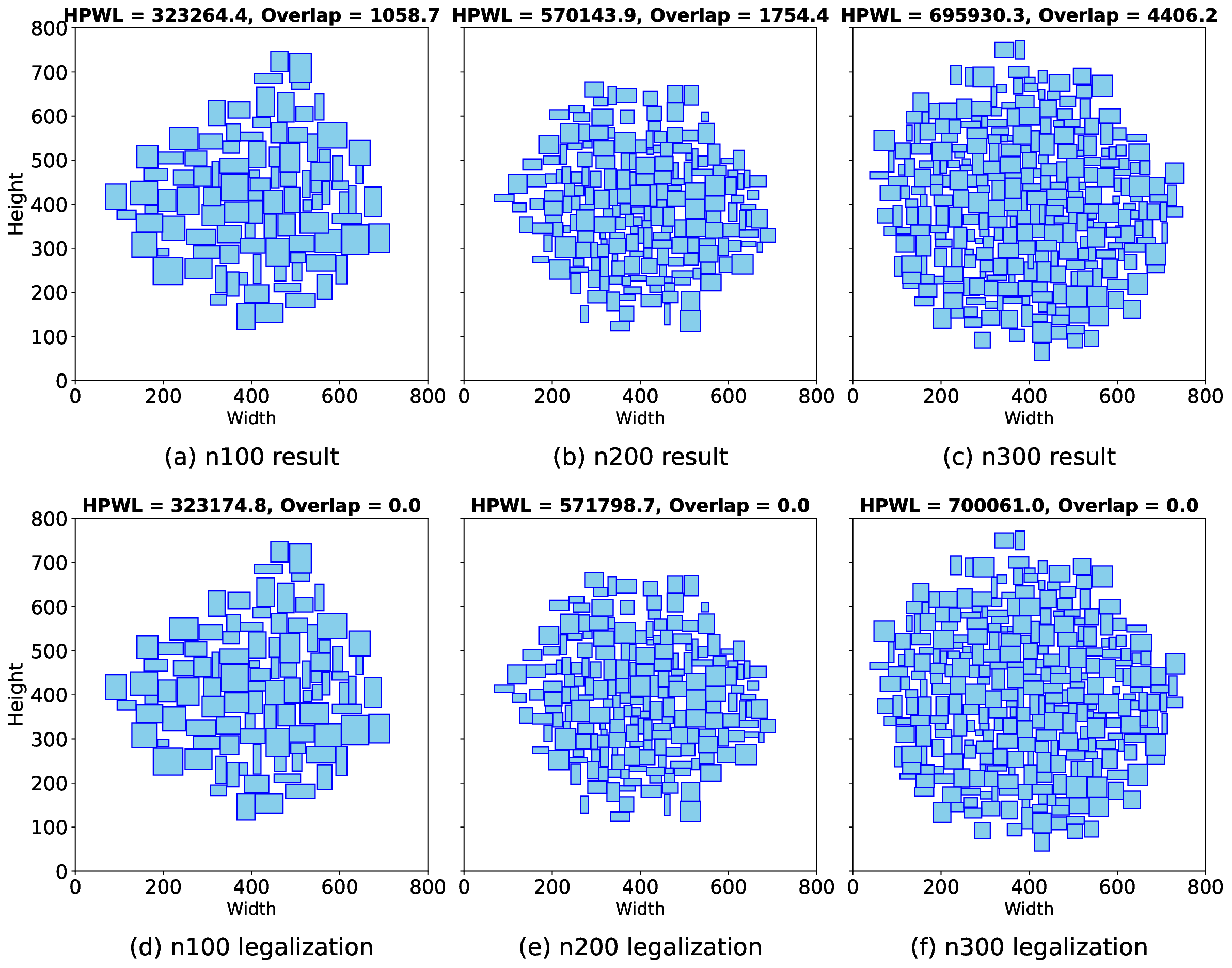}
    \caption{Comparison between the results of RBSM method and the regularization results}
    \label{fig:images2}
\end{figure}

Here we present some figures to display visual placement effects. Fig.~\ref{fig:n10_images2}--Fig.~\ref{fig:n300_images2} show the comparisons of placement effects. Fig.~\ref{fig:images2} shows the legalization results.


\section{Effect of smaller placement regions}

\begin{figure}[ht]
    \centering
    \includegraphics[width=\textwidth]{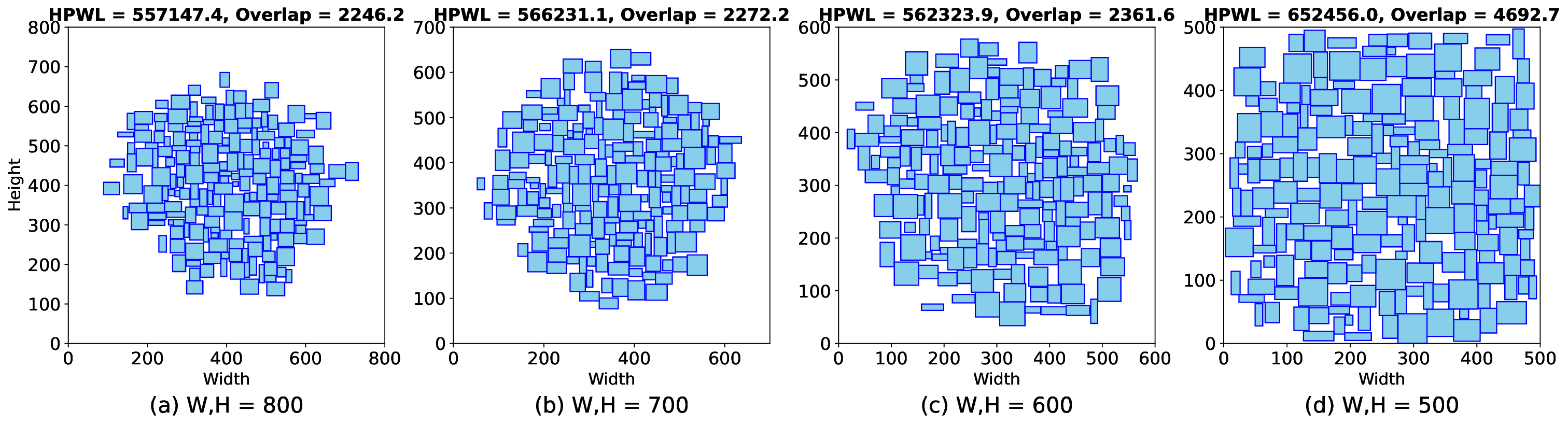}
    \caption{Results of RBSM under different $W,H$ for the $n{=}200$ circuit.}
    \label{fig:images3}
\end{figure}

\begin{figure}[ht]
    \centering
    \includegraphics[width=\textwidth]{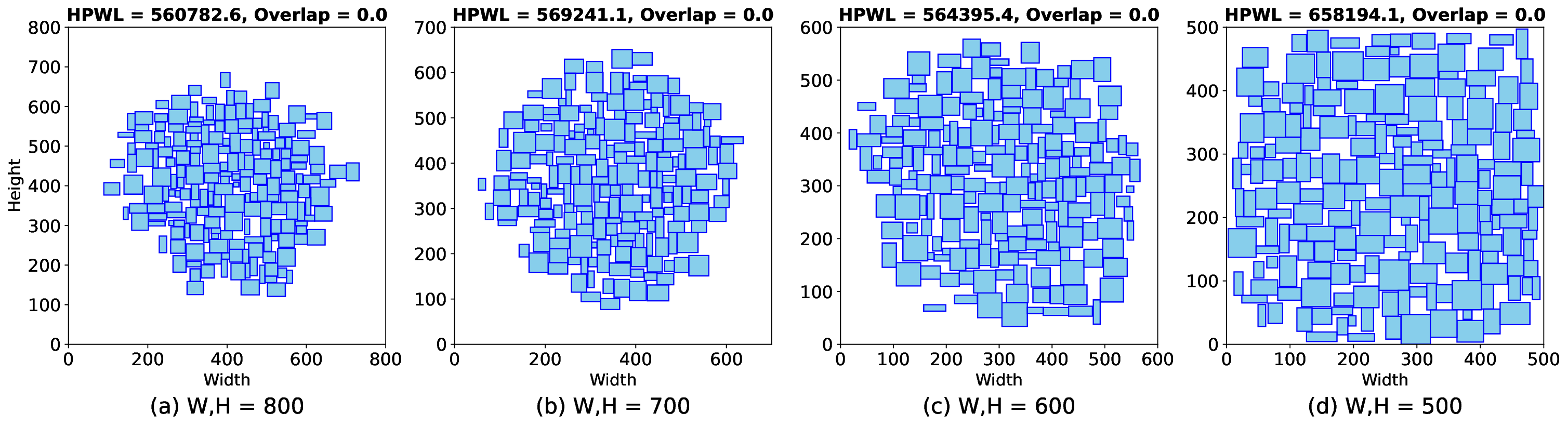}
    \caption{Results after legalization under different $W,H$ for the $n{=}200$ circuit.}
    \label{fig:images4}
\end{figure}

To further evaluate the robustness of our approach, we conducted additional experiments by varying the placement region size $(W,H)$. In particular, we tested smaller region sizes beyond the default $W,H=800$. The results for the $n=200$ circuit are summarized in Fig.~\ref{fig:images3}--Fig.~\ref{fig:images4}.

We observe that within a reasonable range—specifically, when the placement area $(W \times H)$ is at least about 1.8 times larger than the total block area—our algorithm consistently achieves low wirelength (HPWL) and overlap. This indicates that the performance is largely insensitive to the exact choice of $W,H$ under moderate scaling of the placement region. 

However, when the placement region becomes too tight (e.g., $W \times H$ close to the total block area), the wirelength grows significantly, and overlap also increases, since the blocks nearly fill the entire region. This behavior is consistent with practical placement scenarios, where excessively small regions lead to infeasible or poor-quality layouts. 

\end{appendices}   

\bibliography{sn-bibliography}% common bib file

@book{alpert2008handbook,
	title     = {Handbook of Algorithms for Physical Design Automation},
	author    = {Alpert, Charles J. and Mehta, Dinesh P. and Sapatnekar, Sachin S.},
	year      = {2008},
	publisher = {CRC Press},
	address   = {New York},
	isbn      = {978-0-8493-7242-1},
	doi       = {10.1201/9781420013481}
}

@book{kahng2011vlsi,
	title     = {VLSI Physical Design: From Graph Partitioning to Timing Closure},
	author    = {Kahng, Andrew B. and Lienig, Jens and Markov, Igor L. and Hu, Jin},
	year      = {2011},
	publisher = {Springer},
	address   = {Heidelberg, Germany},
    doi       = {10.1007/978-90-481-9591-6}
}

@inproceedings{lasalle2013multi,
	title        = {Multi-threaded graph partitioning},
	author       = {LaSalle, Dominique and Karypis, George},
	booktitle    = {2013 IEEE 27th International Symposium on Parallel and Distributed Processing},
	pages        = {225--236},
	year         = {2013},
	organization = {IEEE},
    doi          = {10.1109/IPDPS.2013.70}
}

@inproceedings{lin2019dreamplace,
	title     = {Dreamplace: Deep learning toolkit-enabled gpu acceleration for modern vlsi placement},
	author    = {Lin, Yibo and Dhar, Shounak and Li, Wuxi and Ren, Haoxing and Khailany, Brucek and Pan, David Z},
	booktitle = {Proceedings of the 56th Annual Design Automation Conference 2019},
	pages     = {1--6},
	year      = {2019},
    doi       = {10.1109/TCAD.2020.3003843}
}

@inproceedings{chan2005multilevel,
	title     = {Multilevel generalized force-directed method for circuit placement},
	author    = {Chan, Tony and Cong, Jason and Sze, Kenton},
	booktitle = {Proceedings of the 2005 international symposium on physical design},
	pages     = {185--192},
	year      = {2005},
    doi       = {10.1109/IPDPS.2013.70}
}

@misc{naylor2001non,
	title     = {Non-linear optimization system and method for wire length and delay optimization for an automatic electric circuit placer},
	author    = {Naylor, William C and Donelly, Ross and Sha, Lu},
	year      = {2001},
	month     = oct # {~9},
	publisher = {Google Patents},
	note      = {US Patent 6,301,693},
    doi       = {US6301693}
}

@article{zhu2015nonsmooth,
	title     = {Nonsmooth optimization method for VLSI global placement},
	author    = {Zhu, Wenxing and Chen, Jianli and Peng, Zheng and Fan, Genghua},
	journal   = {IEEE Trans. Comput.-Aided Des. Integr. Circuits Syst.},
	volume    = {34},
	number    = {4},
	pages     = {642--655},
	year      = {2015},
	publisher = {IEEE},
    doi       = {10.1109/TCAD.2015.2394484}
}

@article{paszke2019pytorch,
	author  = {Paszke, Adam and Gross, Sam and Massa, Francisco and Lerer, Adam and Bradbury, James and Chanan, Gregory and Killeen, Trevor and Lin, Zeming and Gimelshein, Natalia and Antiga, Luca and others},
	title   = {PyTorch: An Imperative Style, High-Performance Deep Learning Library},
	journal = {Adv. Neural Inf. Process. Syst.},
	volume  = {32},
	year    = {2019},
	doi     = {10.48550/arXiv.1912.01703}
}

@article{chang2009essential,
	title     = {Essential issues in analytical placement algorithms},
	author    = {Chang, Yao-Wen and Jiang, Zhe-Wei and Chen, Tung-Chieh},
	journal   = {IPSJ Trans. Syst. LSI Des. Methodol.},
	volume    = {2},
	pages     = {145--166},
	year      = {2009},
	publisher = {Information Processing Society of Japan},
    doi       = {10.2197/ipsjtsldm.2.145}
}

@misc{Gsrc,
	title  = {Gsrc floorplan benchmarks},
	author = {Salama, H. F. and Reeves, D. S. and Viniotis, Y.},
	year   = {1997},
	month  = {June},
}

@book{clarke1990optimization,
	title     = {Optimization and Nonsmooth Analysis},
	author    = {Clarke, Frank H.},
	year      = {1990},
	publisher = {SIAM},
	address   = {Philadelphia, PA}
}

@book{sechen2012vlsi,
	title     = {VLSI Placement and Global Routing Using Simulated Annealing},
	author    = {Sechen, Carl},
	volume    = {54},
	year      = {2012},
	publisher = {Springer Science \& Business Media},
	address   = {New York, NY}
}

@article{huang2021optimization,
	title   = {Optimization models and algorithms for placement of very large scale integrated circuits},
	author  = {Huang, ZP and Li, XQ and Zhu, WX},
	journal = {Operations Research Transactions},
	volume  = {25},
	number  = {03},
	pages   = {15--36},
	year    = {2021},
    doi = {10.15960/j.cnki.issn.1007-6093.2021.03.002}

}

@inproceedings{hsu2011tsv,
	title     = {TSV-aware analytical placement for 3D IC designs},
	author    = {Hsu, Meng-Kai and Chang, Yao-Wen and Balabanov, Valeriy},
	booktitle = {Proceedings of the 48th Design Automation Conference},
	pages     = {664--669},
	year      = {2011}
}

@inproceedings{sun2023floorplanning,
	title        = {Floorplanning of VLSI by Mixed-Variable Optimization},
	author       = {Sun, Jian and Cheng, Huabin and Wu, Jian and Zhu, Zhanyang and Chen, Yu},
	booktitle    = {International Symposium on Intelligence Computation and Applications},
	pages        = {137--151},
	year         = {2023},
	organization = {Springer}
}

@article{chen2008ntuplace3,
	title     = {NTUplace3: An analytical placer for large-scale mixed-size designs with preplaced blocks and density constraints},
	author    = {Chen, Tung-Chieh and Jiang, Zhe-Wei and Hsu, Tien-Chang and Chen, Hsin-Chen and Chang, Yao-Wen},
	journal   = {IEEE Trans. Comput. Aided Des. Integr. Circuits Syst.},
	volume    = {27},
	number    = {7},
	pages     = {1228--1240},
	year      = {2008},
	publisher = {IEEE},
    doi       = {10.1109/TCAD.2008.923063}
}

@inproceedings{caldwell1999optimal,
	title     = {Optimal partitioners and end-case placers for standard-cell layout},
	author    = {Caldwell, Andrew E and Kahng, Andrew B and Markov, Igor L},
	booktitle = {Proceedings of the 1999 international symposium on physical design},
	pages     = {90--96},
	year      = {1999},
    doi       = {10.1109/43.892854}
}

@inproceedings{viswanathan2007fastplace,
	title        = {FastPlace 3.0: A fast multilevel quadratic placement algorithm with placement congestion control},
	author       = {Viswanathan, Natarajan and Pan, Min and Chu, Chris},
	booktitle    = {2007 Asia and South Pacific Design Automation Conference},
	pages        = {135--140},
	year         = {2007},
	organization = {IEEE},
    doi          = {10.1109/ASPDAC.2007.357975}
}

@article{jin2022random,
	author = {Jin, Shi and Li, Lei and Sun, Yiqun},
    title = {On the Random Batch Method for Second Order Interacting Particle Systems},
    journal = {Multiscale Model. Simul.},
    volume = {20},
    number = {2},
    pages = {741-768},
    year = {2022},
    doi = {10.1137/20M1383069}
}

@article{zhou2007bipartite,
	title     = {Bipartite network projection and personal recommendation},
	author    = {Zhou, Tao and Ren, Jie and Medo, Mat{\'u}{\v{s}} and Zhang, Yi-Cheng},
	journal   = {Phys. Rev. E Stat. Nonlin. Soft Matter Phys.},
	volume    = {76},
	number    = {4},
	pages     = {046115},
	year      = {2007},
    doi       = {10.1103/PhysRevE.76.046115},
	publisher = {APS}
}

@inproceedings{lu2014eplace,
	title     = {ePlace: Electrostatics based placement using Nesterov's method},
	author    = {Lu, Jingwei and Chen, Pengwen and Chang, Chin-Chih and Sha, Lu and Huang, Dennis J- H and Teng, Chin-Chi and Cheng, Chung-Kuan},
	booktitle = {Proceedings of the 51st Annual Design Automation Conference},
	pages     = {1--6},
	year      = {2014},
    doi       = {10.1145/2593069.2593216}
}

@article{xiao2021artificial,
	title     = {Artificial bee colony algorithm based on adaptive neighborhood search and Gaussian perturbation},
	author    = {Xiao, Songyi and Wang, Hui and Wang, Wenjun and Huang, Zhikai and Zhou, Xinyu and Xu, Minyang},
	journal   = {Appl. Soft Comput.},
	volume    = {100},
	pages     = {106955},
    issn = {1568-4946},
	year      = {2021},
	publisher = {Elsevier},
	doi       = {https://doi.org/10.1016/j.asoc.2020.106955}
}

@article{weiss1907hypothese,
	title   = {L'hypoth{\`e}se du champ mol{\'e}culaire et la propri{\'e}t{\'e} ferromagn{\'e}tique},
	author  = {Weiss, Pierre},
	journal = {J. Phys. Theor. Appl.},
	volume  = {6},
	number  = {1},
	pages   = {661--690},
	year    = {1907},
    doi     = {10.1051/jphystap:019070060066100}
}

@book{nocedal1999numerical,
	title     = {Numerical Optimization},
	author    = {Nocedal, Jorge and Wright, Stephen J.},
	year      = {1999},
	publisher = {Springer},
	address   = {New York, NY},
    doi       = {10.1007/978-0-387-40065-5}
}

@article{kingma2014adam,
	title = {Adam: a method for stochastic optimization},
	shorttitle = {Adam},
	url = {http://arxiv.org/abs/1412.6980},
	urldate = {2018-04-24},
	journal = {arXiv:1412.6980 [cs]},
	author = {Kingma, Diederik P. and Ba, Jimmy},
	month = dec,
	year = {2014},
}

@inproceedings{gu2020dreamplace,
  title={{DREAMPlace} 3.0: Multi-electrostatics based robust {VLSI} placement with region constraints},
  author={Gu, Jiaqi and Jiang, Zixuan and Lin, Yibo and Pan, David Z},
  booktitle={Proceedings of the 39th International Conference on Computer-Aided Design},
  pages={1--9},
  year={2020}
}

@article{liu2023xplace,
  title={Xplace: An extremely fast and extensible placement framework},
  author={Liu, Lixin and Fu, Bangqi and Lin, Shiju and Liu, Jinwei and Young, Evangeline FY and Wong, Martin DF},
  journal={IEEE Transactions on Computer-Aided Design of Integrated Circuits and Systems},
  volume={43},
  number={6},
  pages={1872--1885},
  year={2023},
  publisher={IEEE}
}

@book{boyd2004convex,
  title={Convex optimization},
  author={Boyd, Stephen P and Vandenberghe, Lieven},
  year={2004},
  address   = {New York},
  publisher={Cambridge University Press}
}

@article{hou2025transplace,
  title={TransPlace: Transferable Circuit Global Placement via Graph Neural Network},
  author={Hou, Yunbo and Ye, Haoran and Yang, Shuwen and Zhang, Yingxue and Xu, Siyuan and Song, Guojie},
  journal={arXiv preprint arXiv:2501.05667},
  year={2025}
}

@article{sun2025automatically,
  title={Automatically discovering heuristics in a complex {SAT} solver with large language models},
  author={Sun, Yiwen and Ye, Furong and Chen, Zhihan and Wei, Ke and Cai, Shaowei},
  journal={arXiv preprint arXiv:2507.22876},
  year={2025}
}

@article{yao2025evolution,
  title={Evolution of Optimization Algorithms for Global Placement via Large Language Models},
  author={Yao, Xufeng and Jiang, Jiaxi and Zhao, Yuxuan and Liao, Peiyu and Lin, Yibo and Yu, Bei},
  journal={arXiv preprint arXiv:2504.17801},
  year={2025}
}

@incollection{robbins1971convergence,
  title={A convergence theorem for non negative almost supermartingales and some applications},
  author={Robbins, Herbert and Siegmund, David},
  booktitle={Optimizing Methods in Statistics},
  pages={233--257},
  year={1971},
  address={New York},
  publisher={Elsevier}
}
%% if required, the content of .bbl file can be included here once bbl is generated
%%\input sn-article.bbl 

\end{document}